\documentclass{article}
\usepackage[a4paper,top=3cm,bottom=3cm,left=2cm,right=2cm]{geometry}
\usepackage{amsmath}
\usepackage{amssymb}
\usepackage{amsthm}
\usepackage{color}
\usepackage{multicol}
\usepackage{graphicx}
\usepackage{algpseudocode} 
\usepackage{algorithm} 
\usepackage{picture}
\usepackage{tikz}
\usepackage{accents}
\usepackage{graphicx}
\usepackage{epstopdf}
\usepackage{subcaption} 
\usepackage[figuresright]{rotating}
\usepackage[none]{hyphenat}
\usepackage{xcolor}

\newcommand{\z}{\phantom{0}}
\renewcommand{\d}{\mathrm{d}}
\newcommand{\vect}[1]{\boldsymbol{#1}}

\usepackage{cases}

\numberwithin{equation}{section}
 
\renewcommand{\d}{\mathrm{d}}
\newcommand{\dt}{\mathrm{dt}}

\renewcommand{\S}{\mathcal{S}}

\newcommand{\vecx}{\boldsymbol{x}}

\begin{document}
\sloppy

\title{Space--time Isogeometric Analysis of cardiac electrophysiology}

\author{Paola Francesca Antonietti$^{1,3}$, Luca Dedè$^1$, Gabriele Loli$^{2}$, Monica Montardini$^{2,3},$ \\ Giancarlo Sangalli$^{2,3}$ and Paolo Tesini$^{5,2,4}$}
\date{}
\maketitle

\begin{center}
\small{$^1$ MOX-Dipartimento di Matematica, Politecnico di Milano\\ Piazza Leonardo da Vinci 32 - 20133, Milano, Italia \\ 
\{paola.antonietti, luca.dede\}@polimi.it} \\
\end{center}
\begin{center}
\small{$^2$ \ Dipartimento di Matematica ``F. Casorati'', Universit\`a di Pavia\\ 
Via A. Ferrata 5 - 27100, Pavia, Italia. \\ 
\{gabriele.loli, monica.montardini, giancarlo.sangalli, paolo.tesini\}@unibs.it} 
\end{center}
\begin{center}
\small{$^3$ \ IMATI-CNR ``Enrico Magenes'',  Pavia, Italy.}
\end{center}
\begin{center}
\small{$^4$ \ Dipartimento di Matematica e Applicazioni, Università degli Studi di Milano-Bicocca, \\ Piazza dell’Ateneo Nuovo, 1 - 20126, Milano, Italia}
\end{center}
\begin{center}
\small{$^5$ \ Dipartimento di Ingegneria Civile, Architettura, Territorio, Ambiente e di Matematica (DICATAM), \\ Università degli Studi di Brescia, via Branze, 43 - 25123, Brescia, Italia}
\end{center}

\begin{abstract}

This work proposes a stabilized 
 space--time method for the monodomain equation coupled with the Rogers--McCulloch ionic model, which is widely used to simulate electrophysiological wave propagation in the cardiac tissue. By extending the Spline Upwind method and exploiting low-rank matrix approximations, as well as preconditioned solvers, we achieve both significant computational efficiency and accuracy.
In particular, we develop a formulation that is both simple and highly effective, designed to minimize spurious oscillations and ensuring computational efficiency. We rigorously validate the method's performance through a series of numerical experiments, showing its robustness and reliability in diverse scenarios.
\end{abstract}

\begin{flushleft}
\end{flushleft}

\section{Introduction}

Isogeometric Analysis (IgA), introduced in \cite{Hughes2005}, is an advanced computational technique that extends the classical finite element method by using spline functions or their generalizations for both representing the computational domain and building the basis functions used in the approximation of partial differential equations.  The smooth spline properties of IgA  offer higher accuracy compared to piecewise polynomials (see, e.g., \cite{Evans2009, Bressan2019}) and compatibility with Computer Aided Design (CAD).

Space–time formulations in Isogeometric Analysis (IgA) present a valuable opportunity to leverage the properties of smooth splines in both space and time, as highlighted in   \cite{Takizawa2014} and \cite{Langer2016}. Notably, reference \cite{Langer2016} focuses on the development of a stabilized IgA framework for the heat equation. In   \cite{Montardini2018} and \cite{Loli2020}, the authors introduce preconditioners and solvers to enhance computational efficiency. Furthermore, in \cite{Saade2021}, a continuous space–time IgA formulation to both linear and non-linear elastodynamics is applied. 

This unified approach allows for a more seamless representation of the problem over time, but it also introduces challenges, particularly when sharp layers or fronts are present in the solution. These sharp layers can cause instabilities, such as spurious oscillations that propagate backward in time, leading to unphysical behaviors in the numerical solutions. To address this issue, the Spline Upwind method (SU) has been proposed in \cite{Loli2023}. This method generalizes classical upwinding techniques, as  SUPG \cite{Brooks1982}, to stabilize space--time IgA formulations, ensuring more accurate  results.

In the context of cardiac electrophysiology, the heart electrical activity is modeled by the bidomain or monodomain equations, which describe the propagation of electrical transmembrane electric potential across the heart tissue. These models are crucial for understanding the heart function, particularly by means of the simulation of the transmembrane potential that triggers muscle contraction, with promising applications in diagnostic and predictive medicine and potentially supporting decision-making procedures in clinics (see, e.g., \cite{Franzone2014, Quarteroni2019}). The bidomain equations divides the cardiac tissue into intracellular and extracellular spaces, governed by reaction--diffusion equations, while the monodomain equations is a parabolic nonlinear PDE coupled to a system of ODEs modelling ionic activity at the cellular level. References \cite{Franzone2014, Franzone2012} offer a comprehensive mathematical derivation of these equations.

An accurate representation of the ionic currents through ion channels is essential for realistic modeling of cardiac excitation. This requires coupling systems of ordinary differential equations with the monodomain equations. However, numerical challenges arise, particularly in resolving sharp and traveling potential fronts, which are crucial for predicting heart activation. 

Various strategies have been developed to address these issues, including modified quadrature rules \cite{Krishnamoorthi2013}, mesh adaptivity \cite{Bendahmane2010, Franzone2006, Southern2012}, and high-order spectral element discretizations \cite{Cantwell2014}. These methods aim to achieve a balance between computational efficiency and solution accuracy.

To address these challenges, IgA has been effectively applied in cardiac electrophysiology, where the use of smooth functions enhances the accuracy of front speed simulations.  In \cite{Patelli2017}, an innovative approach that uses IgA on Non-Uniform Rational B-Spline (NURBS) surfaces has been proposed. The authors show that high-degree basis functions allow for precise front propagation with a reduced number of degrees of freedom. In a similar way, in \cite{Pegolotti2019}, the Roger–McCulloch ionic model \cite{Rogers1994} is employed to address bidomain equations on surfaces, showing that the use of regular B-splines effectively promotes the accuracy of the front propagation velocities and action potential shapes. 
In \cite{torre2022efficient}, the authors propose an isogeometric collocation method to solve the monodomain reaction-diffusion equation, while a multipatch isogeometric method for the monodomain model  has been proposed in \cite{bucelli2021multipatch}.

The aim of this work is to propose a stabilized 
space-time method for the monodomain equation coupled with the Rogers–McCulloch ionic model, by adapting the Spline Upwind method proposed in \cite{Loli2023}, originally developed for the heat equation, to make it computationally efficient for simulating electrophysiological wave propagation in cardiac tissue. 
 The use of low-rank matrix approximation techniques and of efficient preconditioned solvers enable us to tackle challenging problems that would otherwise demand excessive computational resources within the space-time framework.\\
In particular, preconditioners play a crucial role especially in the context of IgA space–time simulations on 3D spatial domains, where the complexity and size of the resulting linear systems can significantly impact on computational performance. Indeed, the application of preconditioners accelerates the iterative solver convergence of resulting linear system. For this reason, we propose a suitable preconditioner, based on the previous construction of \cite{Loli2020} and \cite{Loli2022}.

To evaluate the performance of the method, we present various numerical tests to investigate its numerical stability and computational cost. We also present a comparison with $C^0$ continuity splines and SUPG stabilization. Indeed low continuity splines ensure compatibility with well-established stabilization techniques, such as SUPG, which are extensively documented in the literature and can serve as a useful comparison for our method.

The outline of the paper is as follows.  The basics of IgA are presented in Section \ref{sec:preliminaries}. 
In Section \ref{sec:SU_CE} we apply the stabilization method to the monodomain equation coupled with the Rogers--McCulloch ionic model.
Numerical tests are analyzed in Section \ref{sec:results_ce}, while in Section \ref{sec:conclusion} we draw conclusions.


\section{Preliminary notions of IgA} \label{sec:preliminaries}

We present the notation we use in subsequent sections.

Given $n$ and $p$ two positive integers, we consider the open knot vector 
\begin{equation*}
	\widehat{\Xi}:=\left\{ 0=\widehat{\xi}_1=\dots=\widehat{\xi}_{p+1} \leq \dots \leq \widehat{\xi}_{n}=\dots=\widehat{\xi}_{n+p+1}=1\right\},
\end{equation*}
and the vector $\widehat{Z}:=\left\{\widehat{\zeta}_1, \dots, \widehat{\zeta}_m\right\}$ of knots without repetitions, i.e., breakpoints.

The univariate spline space is defined as
\begin{equation*}
	\widehat{\S}_h^p : = \mathrm{span}\{\widehat{b}_{i,p}\}_{i = 1}^n,
\end{equation*}
where $\widehat{b}_{i,p}$ are the univariate B-splines and $ h$ is the  mesh-size, i.e.,
\begin{equation*}
	h:=\max\{ |\widehat{\xi}_{i+1}-\widehat{\xi}_i| \ \text{s.t.} \ i=1,\dots,n+p \}.
\end{equation*} 

To each univariate B-spline $\widehat{b}_{i,p}$, we associate a Greville abscissa, defined as
\begin{equation*} 
 \gamma_{i,p}:= \frac{\widehat{\xi}_{i+1} + \cdots + \widehat{\xi}_{i+p}}{p} \ \ \ \text{for} \ \ i=1,\dots,n.
\end{equation*} 

Multivariate B-splines are tensor product of univariate B-splines.  

Just for the sake of simplicity,  we consider splines of the same polynomial degree $p_s\geq 1$  in all parametric spatial directions, while  $p_t\geq 1$ is the spline degree in time direction, and $\vect{p}:=\left( p_s,p_t\right)$.

Given positive integers $d$, $n_l$, for $l=1,\dots,d$ and $n_t$, we define $d+1$ univariate knot vectors $\widehat{\Xi}_l:=\left\{ \widehat{\xi}_{l,1} \leq \dots \leq \widehat{\xi}_{l,n_l+p_s+1}\right\}$  for $l=1,\ldots, d$ and $\widehat{\Xi}_t:=\left\{\widehat{\xi}_{t,1} \leq \dots \leq \widehat{\xi}_{t,n_t+p_t}+1\right\}$ and $d+1$ breakpoints vectors $\widehat{Z}_l:=\left\{ \widehat{\zeta}_{l,1},\dots,\widehat{\zeta}_{l,m_l}\right\}$  for $l=1,\ldots, d$ and $\widehat{Z}_t:=\left\{\widehat{\zeta}_{t,1},\dots,\widehat{\zeta}_{t,m_t}\right\}$.
Moreover, let $h_l$ be the mesh-size associated to  $\widehat{\Xi}_l$ for $l=1,\dots,d$, and let $h_t$ be the mesh-size  associated to $\widehat{\Xi}_t$.

The multivariate B-splines are defined as
\begin{equation*} 
	\widehat{B}_{ \vect{i},\vect{p}}(\vect{\eta},\tau) : = \widehat{B}_{\vect{i}_s, p_s}(\vect{\eta}) \widehat{b}_{i_t,p_t}(\tau),
\end{equation*}
where 
\begin{equation*} 
  \widehat{B}_{\vect{i}_s,p_s}(\vect{\eta}):=\widehat{b}_{i_1,p_s}(\eta_1) \ldots \widehat{b}_{i_d,p_s}(\eta_d),
\end{equation*}
$\vect{i}_s:=(i_1,\dots,i_d)$, $\vect{i}:=(\vect{i}_s, i_t)$  and  $\vect{\eta} = (\eta_1, \ldots, \eta_d)$.  
The corresponding spline space is defined as
\begin{equation*}
	\widehat{\S} ^{\vect{p}}_{ {h}  }  := \mathrm{span}\left\{\widehat{B}_{\vect{i}, \vect{p}} \ \middle| \ i_l = 1,\dots, n_l \text{ for } l=1,\dots,d; i_t=1,\dots	,n_t \right\},
\end{equation*} 
where $h:=\max\{h_s,h_t\}$.
We have that
$\widehat{\S} ^{\vect{p}}_{ {h}}  =\widehat{\S} ^{ p_s}_{ {h}_s} \otimes \widehat{\S} ^{p_t}_{h_t}, $ where 
 \[\widehat{\S} ^{p_s}_{h_s} := \mathrm{span}\left\{\widehat{B}_{\vect{i}_s,p_s} \ \middle|  \ i_l =  1,\dots, n_l; l=1,\dots,d  \right\},\]
 is the space of tensor-product splines on $\widehat{\Omega}:=(0,1)^d$, and   
  \[ \widehat{\mathcal{S}}^{ p_t}_{h_t} := \ \text{span}\left\{ \widehat{b}_{i_t,p_t} \ \middle| \ i_t = 1,\dots , n_t\ \right\}.\]

Following \cite{Loli2023}, we assume that $\widehat{\S} ^{p_s}_{h_s}
	\subset C^0(\widehat{\Omega}  )$ and  $\widehat{\S}^{p_t}_{h_t}
	\subset C^{p_t-1}\left((0,1)\right)$. The variable
        continuity in space may be useful for the geometry representation. On the other hand, we consider only maximum continuity with respect to time in order to take advantage of the approximation properties of smooth splines, see, e.g. \cite{Evans2009, Bressan2019}.

We denote by $\Omega\times (0,T)$ the space--time computational domain, where $\Omega\subset\mathbb{R}^d$   and $\Omega$ is parametrized by  $\vect{F}: \widehat{\Omega} \rightarrow {\Omega}$, with  $\vect{F}\in  \left[{\widehat{\mathcal{S}}^{p_s}_{{h}_s}}\right]^d$, and $T>0$ is the final time.
The space--time domain is parametrized by $\vect{G}:\widehat{\Omega}\times(0,1)\rightarrow
 \Omega\times(0,T)$, such that $ \vect{G}(\vect{\eta}, \tau):=(\vect{F}(\vect{\eta}), T\tau )=(\vecx,t).$
 
The spline space with only initial condition, in parametric coordinates, is  
\begin{equation*}
\widehat{\mathcal{X}}_{h}:=\left\{ \widehat{v}_h\in \widehat{\mathcal{S}}^{\vect{p}}_h \ \middle| \  \widehat{v}_h = 0 \text{ on } \widehat{\Omega}\times\{0\} \right\}.
\end{equation*}
We also have that  $\widehat{\mathcal{X}}_{h} =   \widehat{\S} ^{p_s}_{h_s}     \otimes  \widehat{\mathcal{X}}_{t,h_t}  $, where
 \begin{equation*}
   \widehat{\mathcal{X}}_{t,h_t}  := \left\{ \widehat{w}_h\in \widehat{\mathcal{S}}^{ p_t}_{h_t} \ \middle|  \ \widehat{w}_h( 0)=0 \right\}=\text{span}\left\{ \widehat{b}_{i_t,p_t} \ \middle| \ i_t = 2,\dots , n_t\ \right\}.
 \end{equation*}
With colexicographical re-orderings of the basis functions, we write  
\begin{equation}
\widehat{\mathcal{X}}_{h}=\text{span} \left\{ \widehat{B}_{{i}, \vect{p}} \ \middle|\ i=1,\dots,N_{\text{dof}} \right\},
\label{eq:all_basis2}
\end{equation}
 \[\widehat{\S} ^{p_s}_{h_s} := \mathrm{span}\left\{\widehat{B}_{i_s,p_s} \ \middle|  \ i_s =  1,\dots,N_s  \right\},\]
and
  \[ \widehat{\mathcal{S}}^{ p_t}_{h_t} := \ \text{span}\left\{ \widehat{b}_{i_t,p_t} \ \middle| \ i_t = 1,\dots , N_t\ \right\},\]
where $N_{\text{dof}}:=N_{s} N_t$, $N_{s}:=\prod_{l=1}^dn_{l}$, with $l=1,\dots,d$ and $N_t:=n_t-1$.\\
Our isogeometric space is the isoparametric push-forward of \eqref{eq:all_basis2} through the geometric map $\vect{G}$, i.e.,
\begin{equation}
\mathcal{X}_{h} := \text{span}\left\{  B_{i, \vect{p}}:=\widehat{B}_{i, \vect{p}}\circ \vect{G}^{-1} \ \middle| \ i=1,\dots , N_{\text{dof}}   \right\},
\label{eq:disc_space}
\end{equation}
where we can write  
 $\mathcal{X}_{h}=\mathcal{X}_{s,h_s}   \otimes \mathcal{X}_{t,h_t} $,
   with 
\begin{equation*} 
 \mathcal{X}_{s,h_s}    :=\text{span}\left\{ {B}_{i, p_s}:= \widehat{B}_{i, p_s}\circ \vect{F}^{-1} \ \middle| \ i=1,\dots,N_{s} \right\},
\end{equation*}
and
\begin{equation*} 
 \mathcal{X}_{t,h_t}   :=\text{span}\left\{  {b}_{i,p_t}(\cdot):= \widehat{b}_{i,p_t}\left(\frac{\cdot}{T}\right) \ \middle| \ i=1,\dots,N_t \right\}.
\end{equation*}

Finally, following \cite{Daveiga2014}, we define the support extension for $\vect{i}_s:=(i_1,\dots,i_d)$, with $i_l=1,\dots,n_l$, $l=1,\dots,d$ and $i_t=1,\dots,N_t$, as
\begin{equation} \label{suppext}
	\widetilde{I}_{\vect{i}_s,i_t}:=\left( \widehat{\xi}_{1,i_1-p_s} \ \text{,} \ \widehat{\xi}_{1,i_1+p_s+1}\right) \times \cdots  \times \left( \widehat{\xi}_{d,i_d-p_s} \ \text{,} \ \widehat{\xi}_{d,i_d+p_s+1}\right) \times \left( \widehat{\xi}_{t,i_t-p_t} \ \text{,} \ \widehat{\xi}_{t,i_t+p_t+1}\right).
\end{equation}

To further explore the properties of B-splines and their application in IgA, we refer the reader to \cite{Cottrell2009}.


\section{Spline Upwind for cardiac electrophysiology} \label{sec:SU_CE}

In the context of cardiac electrophysiology, we consider the monodomain equation with the Rogers–McCulloch ionic model, proposed in \cite{Rogers1994}. The dimensionless unknowns are the transmembrane potential $u$ and the recovery variable $w$, and the governing differential problem is
\begin{equation}
\label{eq:McCulloch}
	\left\{
	\begin{array}{rcllrcl}
		 C_m \partial_t u -  \nabla\cdot(\mathbf{D} \nabla u)  + c_1u(u-a)(u-1) + c_2 u w & = & f & \mbox{in } &\Omega \!\!\!\! &\times & \!\!\!\! (0, T), \\[1pt]
		 \partial_t w -  b(u-d_ew) & = & 0 & \mbox{in } &\Omega \!\!\!\! &\times & \!\!\!\! (0, T), \\[1pt]
		 \vect{n}\cdot(\mathbf{D} \nabla u)   & = & 0 & \mbox{on } &\partial\Omega \!\!\!\! &\times& \!\!\!\! [0, T], \\[1pt]	 	
		 u & = & 0 & \mbox{in } &\Omega\!\!\!\!  &\times &\!\!\!\! \lbrace 0 \rbrace ,\\[1pt]	
		 w & = & 0 & \mbox{in } &\Omega\!\!\!\!  &\times &\!\!\!\! \lbrace 0 \rbrace, 
	\end{array}
	\right.
\end{equation}
where $\vect{n}$ is the exterior normal, $a$, $b$, $c_1$, $c_2$ and $d_e$ are dimensionless positive parameters, specific of the Rogers–McCulloch model.
The dimensionless constant  $C_m$ is the local membrane capacitance, while $\mathbf{D}=\mathbf{D}(\vect{x})\in L^{\infty}(\overline{\Omega};\mathbf{R}^{d\times d})$ is the conductivity tensor, that we assume to be symmetric and uniformly positive definite. 
We assume that  $f\in L^2(\Omega  \times  (0, T))$.

We consider the standard Galerkin method\\
find $u_h \text{,} w_h \in \mathcal{X}_{h}$ such that
\begin{equation}
\label{eq:McCulloch_gal}
	\left\{
	\begin{array}{rcllrcl}
		 \mathcal{A}(u_h,w_h;v_h) & = & \mathcal{F}(v_h) & \ \ \forall v_h  \in  \mathcal{X}_{h},  \\[1pt]
		 \mathcal{L}(w_h,u_h;z_h) & = & 0 & \ \ \forall z_h  \in  \mathcal{X}_{h}, \\
	\end{array}
	\right.
\end{equation}
where
\begin{equation*}
 \mathcal{A}(u_h,w_h;v_h) := \int_{0}^T\int_{\Omega} \big{(}C_m \partial_t {u_h} \, v_h +    (\nabla u_h)^T\mathbf{D}   \nabla v_h  + \left(c_1(u_h-a)(u_h-1) + c_2 w_h \right)\, u_h\, v_h  \big{)}\,\d\Omega\,   \dt ,
\end{equation*}
\begin{equation}
\mathcal{F}(v_h)  := \int_0^T\int_{\Omega}f\,   v_h \,\d\Omega\,\dt ,
\label{eq:rhsfun}
\end{equation}
and
\begin{equation*}
 \mathcal{L}(w_h,u_h;z_h) := \int_{0}^T\int_{\Omega}  \big{(} \partial_t {w_h} \, z_h -  b(u_h-d_e w_h) \, z_h \big{)}\,\d\Omega\,   \dt .
\end{equation*}

We also consider the following stabilized formulation, based on the Spline Upwind technique introduced in \cite{Loli2023}: find $u_h, w_h \in \mathcal{X}_{h}$ such that
\begin{equation}
	\label{eq:McCulloch_SU}
	\left\{
	\begin{array}{rcllrcl}
		\mathcal{A}(u_h,w_h;v_h)+\mathcal{S}(u_h ; v_h) & = & \mathcal{F}(v_h) & \ \ \forall v_h  \in  \mathcal{X}_{h},  \\[1pt]
		\mathcal{L}(w_h,u_h;z_h) & = & 0 & \ \ \forall z_h  \in  \mathcal{X}_{h}, \\
	\end{array}
	\right.
\end{equation}
where the stabilization term $\mathcal{S}(u_h ; v_h)$ corresponds to $\mathcal{S}_{\mathrm{SU,2}}(u_h ; v_h)$ introduced in~\cite{Loli2023}, but with a different $\theta$. In particular, $\mathcal{S}(u_h ; v_h)$ is given by:
\begin{equation*}\label{eq:SU_CE_u}
	\mathcal{S}(u_h, w_h ; v_h) := C_m \sum_{k=1}^{p}{\int_{0}^{T} \tau_{k}(t) \int_{\Omega} {\theta (u_h,w_h) \,\partial_t^k u_h\partial_t^k v_h} \ \d\Omega \ \dt},
\end{equation*}
where, for $k=1,\ldots,p_t$, $\tau_{k}(T \ \cdot) \in \widehat{\mathcal{S}}^{ p_t-k}_{h_t}\subset C^{p_t-k-1}$, and are chosen such that:
\begin{equation}\label{tau}
	\int_0^T{ b'_{\ell+i,p_t} b_{i,p_t} \ \dt} + \sum_{k=1}^{p_t}\int_{0}^{T}{\tau_{k}(t)\,{ b^{(k)}_{\ell+i,p_t} b^{(k)}_{i,p_t}} \ \dt} = 0,
\end{equation}
for $i = 1,\ldots,N_t-1,\ \ell=1,\ldots,r,$ with $r=\min(p_t,N_t-i),$ as in \cite{Loli2023}. Additionally, $\theta(u_h,v_h)\equiv\theta(\vect{x},t,u_h,w_h)$ denotes the function defined as the $(d+1)$-linear interpolation computed at the Greville abscissae mapped to the physical domain through the map $\mathbf{G}$ of the following values:
\begin{equation}\label{eq:theta_ce_u}
	\Theta_{\vect{i}_s,i_t}:=\min \left(\frac{\left\| C_m \partial_t {u_{h}} -   \nabla \cdot (\mathbf{D}\nabla u_{h}) + c_1u_h(u_h-a)(u_h-1) + c_2 u_h w_h - f  \right\|_{L^{\infty}(\widetilde{I}_{\vect{i}_s,i_t})}}{C_m(T^{-1}\left\| {u_{h}} \right\|_{L^{\infty}(\Omega \times[0,T])} + \left\| \partial_t {u_{h}} \right\|_{L^{\infty}(\Omega \times[0,T])})},1\right),
\end{equation}
for \(i_l=1,\ldots,n_l,\ l=1,\ldots,d,\ i_t=1,\ldots,N_t\). 

Unlike in \cite{Loli2023}, the $L^{\infty}$ norms of the residuals in \eqref{eq:theta_ce_u} are computed on the support extensions of the B-splines (see \eqref{suppext}), rather than on two mesh elements. This modification extends the stabilization to a larger portion of the domain, enhancing stability. To further simplify, we just drop the first order (SUPG) consistent stabilization term  from \cite{Loli2023}. Near the layers, high-order stabilization is applied primarily to ensure the time derivative matrix block is (almost) triangular, see  \cite{Loli2023} and below.

\subsection{Non-linear solver}\label{sec:semiimplicit}

We solve \eqref{eq:McCulloch_gal} and \eqref{eq:McCulloch_SU} by the following fixed point scheme:\\
given $u^{0}_h, w^{0}_h \in  \mathcal{X}_{h}$, the $(k+1)$-th iteration consists in finding $\widetilde{u}^{k+1}_h$, $\widetilde{w}^{k+1}_h \in \mathcal{X}_{h}$ such that
\begin{equation} 
	\left\{
	\begin{array}{rcllrcl}
		 \mathcal{A}(\widetilde{u}^{k+1}_h,u_h^{k},w_h^{k};v_h)+ {\mathcal{S}}(\widetilde{u}^{k+1}_h, u_h^k, w_h ; v_h) & = & \mathcal{F}(v_h) & \ \ \forall v_h  \in  \mathcal{X}_{h},  \\[1pt] 
		 \mathcal{L}(\widetilde{w}^{k+1}_h;z_h) & = & \mathcal{G}(u^k_h;z_h) & \ \ \forall z_h  \in  \mathcal{X}_{h}, \\
	\end{array}
	\right. 
	\label{eq:solverSU} 
\end{equation}
where, with abuse of notation, in \eqref{eq:solverSU} we defined
\begin{equation*}
 \mathcal{A}(\widetilde{u}^{k+1}_h,u_h^{k},w_h^{k};v_h)  :=  \int_{0}^T\int_{\Omega}  \left( C_m \partial_t {\widetilde{u}^{k+1}_h} \, v_h  +   (\nabla\widetilde{u}^{k+1}_h)^T\mathbf{D}  \nabla v_h  + C_r(u^k_{h},w^k_h)\,\widetilde{u}^{k+1}_h\, v_h  \right)\,\d\Omega\,   \dt ,
\end{equation*}
with 
\begin{equation}
C_r(u^k_{h},w^k_h):=c_1(u^{k}_h-a)(u^{k}_h-1) + c_2 w^{k}_h,
\label{eq:reaction}
\end{equation}
and
\begin{equation*}
 \mathcal{S}(\widetilde{u}^{k+1}_h,u_h^{k},w_h^{k};v_h)  :=  C_m \sum_{l=1}^{p}{\int_{0}^{T} \tau_{k}(t) \int_{\Omega} {\theta (u^k_h,w^k_h) \,\partial_t^l \tilde{u}^{k+1}\partial_t^l v_h} \ \d\Omega \ \dt},
\end{equation*}
while
\begin{equation*}
 \mathcal{L}(\widetilde{w}^{k+1}_h;z_h) := \int_{0}^T\int_{\Omega}  \left( \partial_t {\widetilde{w}^{k+1}_h} \, z_h +b \, d_e \, \widetilde{w}^{k+1}_h \,z_h \right)\,\d\Omega\,   \dt\quad \text{and}   \quad
 \mathcal{G}(u^{k}_h;z_h) := \int_{0}^T\int_{\Omega}   b \, u^{k}_h    \, z_h\,\d\Omega\,   \dt  ,
\end{equation*} 
and $\mathcal{F}$ is defined in \eqref{eq:rhsfun}.
System \eqref{eq:solverSU} consists of two decoupled equations, each of which can be solved independently.
In order to speed up the computation of the stabilizing term, we consider a low-rank approximation of $\theta(u_h^k,w_h^k)$.   
Given a relative tolerance $\varepsilon>0$, using the algorithm in \cite{Wenjian2018}, we find $R\in \mathbb{N}$ s.t. $0<R\leq \text{min}(N_t,N_s)$, $\mathbf{U}\in \mathbb{R}^{N_t\times R}$, $\mathbf{V}\in \mathbb{R}^{N_s\times R}$ and $\mathbf{R}\in\mathbb{R}^{R\times R}$ a diagonal matrix s.t. 
\begin{equation}\label{eq:tol}
	\frac{\left\| \mathbf{\widehat{\Theta}}-\mathbf{U}\mathbf{R}\mathbf{V}^T \right\|_\text{F}}{\left\| \mathbf{\widehat{\Theta}}\right\|_{\text{F}}} \leq \varepsilon,
\end{equation}
where $\mathbf{\widehat{\Theta}}\in \mathbb{R}^{N_t\times N_s}$ is obtained by reshaping the tensor \eqref{eq:theta_ce_u} and $\left\| \cdot \right\|_{\text{F}}$ is the Frobenius norm. 

In this way, we can write the function $\theta(u_h^k,w_h^k)$ as
\begin{equation*}
\theta(u_h^k,w_h^k)\equiv\theta(u_h^k,w_h^k,\vect{x},t)\approx  \sum_{r=1}^R [ \mathbf{R}]_{r,r} \theta_{t,r}(t)\theta_{s,r}(\vect x),
\end{equation*}
where
$\theta_{t,r}(t)$ is a linear interpolation of the $r$-th column of $\mathbf{U}$, on the Greville abscissae in the time direction and
$\theta_{s,r}(\vect x)$ is a $d$-linear interpolation of $r$-th column of $\mathbf{V}$, on the Greville abscissae in the space directions. Of course, both $\theta_{t,r}$ and $\theta_{s,r}$ have to be computed at each non-linear iteration, as they depend on $u_h^k$ and $w_h^k$, even if, for the sake of simplicity in the notation, this dependence is omitted.

Therefore, we have the following approximation
\begin{align}
 \mathcal{S}(\widetilde{u}^{k+1}_h,u_h^{k},w_h^{k};v_h)&\approx \widetilde{\mathcal{S}}(\widetilde{u}^{k+1}_h,u_h^{k},w_h^{k};v_h) \\ &:= C_m \sum_{r=1}^R [ \mathbf{R}]_{r,r} \sum_{l=1}^{p}  {\int_{0}^{T} \tau_{k}(t)\,  \theta_{t,r}(t) \int_{\Omega} \theta_{s,r}(\vect x)\partial_t^l \tilde{u}^{k+1}_h\partial_t^l v_h} \ \d\Omega \, \dt,
 \label{eq:tildeS}
\end{align}
and the fixed point scheme to be solved becomes: given $u^{0}_h, w^{0}_h \in  \mathcal{X}_{h}$, the $(k+1)$-th iteration consists in finding $\widetilde{u}^{k+1}_h$, $\widetilde{w}^{k+1}_h \in \mathcal{X}_{h}$ such that
\begin{equation} 
	\left\{
	\begin{array}{rcllrcl}
		 \mathcal{A}(\widetilde{u}^{k+1}_h,u_h^{k},w_h^{k};v_h)+ \widetilde{\mathcal{S}}(\widetilde{u}^{k+1}_h,u_h^{k},w_h^{k};v_h) & = & \mathcal{F}(v_h) & \ \ \forall v_h  \in  \mathcal{X}_{h},  \\[1pt] 
		 \mathcal{L}(\widetilde{w}^{k+1}_h;z_h) & = & \mathcal{G}(u^k_h;z_h) & \ \ \forall z_h  \in  \mathcal{X}_{h}. \\
	\end{array}
	\right. 
	\label{eq:solverSU_approx} 
\end{equation}

Moreover, in order to promote the convergence of fixed point iterations, we introduce relaxation 

\begin{equation*}
	u^{k+1}_h = \alpha \widetilde{u}_h^{k+1} +(1-\alpha) u_h^{k} \ \ \ \text{and} \ \ \ w_h^{k+1} = \alpha \widetilde{w}_h^{k+1} +(1-\alpha) w_h^{k},
\end{equation*}
where we set $\alpha=0.5$. 

Finally, we remark that the functions \(\tau_k\) appearing in \eqref{tau} can be computed only once and, since \(C_m\) is a positive constant, the high-order stabilization ensures that the   time derivative matrix block is (almost) a lower-triangular matrix at each fixed-point iteration.

\subsubsection{Discrete linear system}
The linear systems resulting from \eqref{eq:solverSU_approx} are
\begin{align}
\mathbf{A}\widetilde{\mathbf{u}}^{k+1}&=\mathbf{f}, \label{eq:lin_sys_A}\\
\mathbf{L}\widetilde{\mathbf{w}}^{k+1}&=\mathbf{g}, \label{eq:lin_sys_L}
\end{align}
where 
$
[\mathbf{A}]_{i,j}:=\mathcal{A}(B_{j,\vect{p}},u_h^{k},w_h^{k};B_{i,\vect{p}}) , \quad [\mathbf{L}]_{i,j}:=(B_{j,\vect{p}},B_{i,\vect{p}}), \quad [\mathbf{f}]_i:=\mathcal{F}(B_{i,\vect{p}}),\quad [\mathbf{l}]_i:=\mathcal{G}(B_{i,\vect{p}}).
$
Exploiting the tensor-product structure of the isogeometric space $\mathcal{X}_h$  \eqref{eq:disc_space}, we have
\begin{align*}
\mathbf{A} &:= C_m \, \mathbf{W}_t \otimes  \mathbf{M}_s + \, \mathbf{M}_t\otimes \mathbf{K}_s + \mathbf{M}_{\text{R}} + C_m \sum_{r=1}^R [ \mathbf{R}]_{r,r} \sum_{l=1}^{p} \mathbf{S}^t_{r,l} \otimes \mathbf{S}^s_{r,l},\\
\mathbf{L} &:= (\mathbf{W}_t +  b \, d_e \, \mathbf{M}_t) \otimes \mathbf{M}_s ,
\end{align*}
where, using the definition \eqref{eq:reaction}, for $i,j=1, \dots, N_{\text{dof}}$  we have
\begin{equation*}
[ \mathbf{M}_{\text{R}}]_{i,j}  =  \int_0^T  \int_{\Omega}  C_r(u^k_{h},w^k_h) \ B_{i, \vect{p}}(\vect{x},t) \  B_{j, \vect{p}}(\vect{x},t) \ \d\Omega \, \dt, 
\end{equation*}
for $i,j=1, \dots, N_t$,
\begin{equation*}
	{[\mathbf{S}^t_{r,l}]}_{i,j} =   {\int_{0}^{T} \tau_{k}(t) \, \theta_{t,r}(t) \, \partial_t^l b_{i, p_t}(t)\, \partial_t^l b_{j, p_t}} \  \dt,
\end{equation*}
\[
 [ \mathbf{W}_t]_{i,j}  = \int_{0}^T   b'_{j,  p_t}(t)\,  b_{i, p_t}(t) \, \dt \quad \text{and} \quad   [ \mathbf{M}_t]_{i,j} = \int_{0}^T\,  b_{i, p_t}(t)\,  b_{j, p_t}(t)  \, \dt , 
\]
while, for $i,j=1, \dots, N_s$,
\begin{equation*}
{[\mathbf{S}^s_{r,l}]}_{i,j}=\int_{\Omega} \theta_{s,r}(\vect x) \, B_{i, p_s}(\vect{x}) B_{j, p_s}(\vect{x}) \d\Omega,
\end{equation*}
\[
[ \mathbf{K}_s]_{i,j}  =  \int_{\Omega} (\nabla B_{i, p_s}(\vect{x}))^T\mathbf{D}(\vect{x}) \nabla B_{j, p_s}(\vect{x}) \quad \text{and} \quad [ \mathbf{M}_s]_{i,j}  =  \int_{\Omega}  B_{i, p_s}(\vect{x}) \  B_{j, p_s}(\vect{x}) \ \d\Omega.
\]

\subsection{Preconditioner}\label{sec:prec}
 We consider for the linear system \eqref{eq:lin_sys_A} the following preconditioner
 $$
[\widehat{\mathbf{A}}]_{i,j}=\widehat{\mathcal{A}}(\widehat{B}_{j,\vect{p}}, \widehat{B}_{i,\vect{p}}) 
 $$
where 
$$
\widehat{\mathcal{A}}(u,v)  :=  \int_{0}^T\int_{\widehat{\Omega}}  \left( C_m \partial_t {{u}} \, v  +   \nabla {u} \cdot \nabla v  + a \, c_1 \, u \, v  \right)\,\d\Omega\,   \dt .
$$
We have that

\begin{equation}
	\widehat{\mathbf{A}}=  C_m \, \mathbf{W}_t \otimes  \widehat{\mathbf{M}}_s + \mathbf{M}_t\otimes(   \, \widehat{\mathbf{K}}_s + a \, c_1 \, \widehat{\mathbf{M}}_s),
	\label{eq:preconditioner}
\end{equation}
where for $i,j=1,\dots,N_s$  
\[
[ \widehat{\mathbf{K}}_s]_{i,j}  =  \int_{\widehat{\Omega}} \nabla \widehat{B}_{i, p_s}(\vect{x}) \cdot \nabla \widehat{B}_{j, p_s}(\vect{x}) \quad \text{and} \quad [ \widehat{\mathbf{M}}_s]_{i,j}  =  \int_{\widehat{\Omega}}  \widehat{B}_{i, p_s}(\vect{x}) \  \widehat{B}_{j, p_s}(\vect{x}) \ \d\widehat{\Omega}
\]
are the stiffness and mass matrices in the parametric domain, respectively. It also holds
$$
\widehat{\mathbf{K}}_s=\sum_{i=1}^d\widehat{M}_d\otimes\dots \widehat{M}_{i+1}\otimes\widehat{K}_i\otimes\widehat{M}_{i}\otimes\dots \otimes \widehat{M}_1\quad\text{and}\quad \widehat{\mathbf{M}}_s=\widehat{M}_d\otimes\dots\otimes \widehat{M}_1,
$$
where $\widehat{K}_l$ and $\widehat{M}_l$ are the univariate stiffness and mass matrices, respectively.
To apply the preconditioner \eqref{eq:preconditioner}, we generalize the technique proposed in \cite{Loli2020}. Thus, following \cite{Loli2020}, we consider the generalized eigendecomposition of the pencils $ (\widehat{\mathbf{K}}_l, \widehat{\mathbf{M}}_l )$ for $l=1,\dots,d$, which gives  the matrices $\mathbf{U}_l$   for $l=1,\dots,d$ for $l=1,\dots,d$, that  contain in each column the $\widehat{\mathbf{M}}_l$-orthonormal generalized eigenvectors   and  $\mathbf{\Lambda}_l$ that are diagonal matrices whose entries contain the corresponding generalized eigenvalues. Moreover, we   define
\[ 
[\mathbf{w}]_i= [\mathbf{W}_t]_{i,N_t} \quad\text{and}  \quad [\mathbf{m}]_i= [\mathbf{M}_t]_{i,N_t}   \quad\text{for}\quad  i=1,\dots, N_t-1, 
\]
\[ 
\z[\overset{\circ}{\mathbf{W}}_t]_{i,j} = [\mathbf{W}_t]_{i,j}    \quad\text{and}  \quad [\overset{\circ}{\mathbf{M}}_t]_{i,j} = [\mathbf{M}_t]_{i,j}   \quad \text{for}   \quad i,j=1,\dots, N_t-1.
\]
and we  consider the matrices $\overset{\circ}{\mathbf{U}}_t$ and $\boldsymbol{\Lambda}_t$, that are the matrix whose columns contain  the $\overset{\circ}{ {\mathbf{M}}_t}$-orthogonal generalized eigenvectors of
 the pencil 
 $(\overset{\circ}{\mathbf{W}}_t,\overset{\circ}{\mathbf{M}}_t)$ and the matrix of the corresponding eigenvalues, respectively.
We then define the matrix $\mathbf{U}_t$   as
\begin{equation*}
   \mathbf{U}_t :=
\begin{bmatrix}
\overset{\circ}{\mathbf{U}}_t &  \mathbf{t} \\[2pt]
   \mathbf{0}^T & \rho \\
\end{bmatrix},
\end{equation*} 
where $\mathbf{0} \in \mathbb{R}^{N_t-1}$ denotes the null vector, while  
\[
\begin{bmatrix}
\mathbf{t}\\ \rho
\end{bmatrix}:=\frac{\begin{bmatrix}
\mathbf{v}\\ 1
\end{bmatrix}}{  \left( [\mathbf{v}^*\ 1]\mathbf{M}_t\begin{bmatrix}
\mathbf{v}\\ 1
\end{bmatrix}\right)^{\tfrac12}  }, 
\]
and $\mathbf{v}\in \mathbb{C}^{N_t-1}$ such that
\[
 \overset{\circ}{\mathbf{M}}_t\mathbf{v}=-\mathbf{m}.
\] 
   
Finally, we set $\mathbf{\Delta}_t:=\mathbf{U}^*_t\mathbf{W}_t\mathbf{U}_t$. The matrix $\mathbf{\Delta}_t$ has an arrowhead structure  (see \cite{Loli2020} for details).  
 
 Then, factorizing the common terms, we get 
that $\widehat{\mathbf{A}}$ can be written as
\begin{equation}
	\widehat{\mathbf{A}}=\left( \mathbf{U}_t^* \otimes \mathbf{U}_s^T \right)^{-1} \left( C_m \mathbf{\Delta}_t  \otimes \mathbb{I}_{N_s} +   \ \mathbb{I}_{N_t} \otimes \mathbf{\Lambda}_s +a \, c_1 \, \mathbb{I}_{N_t}  \otimes \mathbb{I}_{N_s}\right)\left( \mathbf{U}_t \otimes \mathbf{U}_s \right)^{-1},
	\label{eq:factorized}
\end{equation}
 where   $\mathbf{\Lambda}_s:=\sum_{l=1}^d \mathbb{I}_{n_{d}}\otimes\dots\otimes\mathbb{I}_{n_{l+1}} \otimes  \mathbf{\Lambda}_l \otimes \mathbb{I}_{n_{l-1}}\otimes\dots\otimes\mathbb{I}_{n_{1}}$, $\mathbf{U}_s:=\mathbf{U}_d \otimes \cdots \otimes \mathbf{U}_1$ and $\mathbb{I}_n$ denotes the identity of dimension $n\times n$. Note that the second term of \eqref{eq:factorized} exhibits a block-arrowhead structure: 
\begin{equation*}
  \begin{aligned}
   C_m \mathbf{\Delta}_t  \otimes \mathbb{I}_{N_s} +   \ \mathbb{I}_{N_t}  \otimes \mathbf{\Lambda}_s & +a \, c_1 \, \mathbb{I}_{N_t}  \otimes \mathbb{I}_{N_s} =
\begin{bmatrix}
\mathbf{H}_1 &  &  & \mathbf{B}_1 \\[2pt]
   &\ddots & &\vdots \\[2pt]
  \quad &  & \mathbf{H}_{N_t-1} &   \mathbf{B}_{N_t-1} \\[2pt]
-\mathbf{B}^*_1 &  \ldots & -\mathbf{B}^*_{N_t-1} & \mathbf{H}_{N_t}
\end{bmatrix},
  \end{aligned}
\label{eq:arrow}
\end{equation*} 
where $\mathbf{H}_i$ and $\mathbf{B}_i$ are diagonal matrices defined as 
\begin{equation*}
\label{eq:arrow_expl}
\mathbf{H}_i:={  C_m}[\mathbf{\Lambda}_t]_{i,i}\mathbb{I}_{N_s}+\mathbf{\Lambda}_s +a \, c_1 \, \mathbb{I}_{N_s} \quad \text{ and }
 \quad \mathbf{B}_i:={  C_m}[\mathbf{l}]_i\mathbb{I}_{N_s}  \quad  
  \text{ for }\quad i=1,\dots,N_t-1,
\end{equation*}
\begin{equation*} 
 \mathbf{H}_{N_t}:={  C_m}\sigma{\mathbb{I}}_{N_s}+ \mathbf{\Lambda}_s + a \, c_1 \, \mathbb{I}_{N_s},
\end{equation*}  
while $\mathbf{l}:=\overset{\circ}{{\mathbf{U}}}_t^*\left[\overset{\circ}{\mathbf{W}}_t \ \mathbf{w}\right] \begin{bmatrix}\mathbf{t}\\ \rho \end{bmatrix} $ and
$
\sigma:=\left[ {\mathbf{t}}^*   {\rho}^*\right]  \mathbf{W}_t 
\begin{bmatrix}
\mathbf{t}\\ \rho
\end{bmatrix}$.
This property is used in the  application of the preconditioner \eqref{eq:factorized}, that we perform with  \cite[Algorithm 1 Extended FD]{Loli2020}. 
Finally, parallelization presents a further opportunity for enhancing computational efficiency, but is not adopted in this paper.

We now focus on the solution of \eqref{eq:lin_sys_L}: we note that, exploiting the properties of the Kronecker product, the solution can be obtained as
\begin{equation*}
\mathbf{w} = \text{vec}\left( \mathbf{M}_s^{-1} \mathbf{\widetilde{G}} \left( \mathbf{W}_t + b \, d_{e} \mathbf{M}_t \right)^{- \top} \right),
\end{equation*}
where the ``vec'' operator applied to a matrix stacks its columns into a vector, and \(\mathbf{\widetilde{G}}\) is the \(N_s \times N_t\) matrix such that \(\mathbf{g} = \text{vec}(\mathbf{\widetilde{G}})\).
Therefore, \(\mathbf{w}\) can be computed by solving \(N_s\) independent systems associated with the \(N_t \times N_t\) matrix \(\left( \mathbf{W}_t + b\,d_{e} \mathbf{M}_t \right)^{\top}\), as well as \(N_t\) independent systems associated with the \(N_s \times N_s\) spatial mass matrix \(\mathbf{M}_s\). For the latter we can efficiently leverage the preconditioner proposed in \cite{Loli2022}.


\section{Numerical Results} \label{sec:results_ce}

All numerical tests are performed using {\sffamily MATLAB} with {\sffamily GeoPDEs} \cite{Vazquez2016} and {\sffamily Tensorlab} \cite{Sorber2014} on Intel(R) Xeon(R) Gold 6130 CPU processor, running at 2.10 GHz and with 128 GB of RAM.  We use B-splines with the maximum continuity allowed, except where otherwise specified, and  we consider the same polynomial degree for space and time, i.e., $p_s=p_t=:p$.
Following \cite{Rogers1994}, we set $a=0.13$, $b=0.013$, $c_1=0.26$, $c_2=0.1$ and $d_e=1$. 

For the low-rank approximation \eqref{eq:tol}, we use  {\ttfamily svdsketch} {\sffamily MATLAB} function with $\varepsilon=10^{-1}$ and the fixed point method, presented in Section \ref{sec:semiimplicit}, is used to solve nonlinearities in the equations, with the stopping criterion: $\left\| \textbf{u}^{k+1} - \textbf{u}^{k}\right\|_{\infty}\leq \delta$, where $\delta=10^{-4}$ and  $\textbf{u}^n$ is the vector of the  degrees of freedom of the $n$-th iterate of the fixed point scheme for $n\geq 1$.

In the numerical experiments,  $\chi_{[\overline{\psi}_1,\overline{\psi}_2]}(\psi)$ refers to the characteristic function, defined as 
\begin{equation*} 
\chi_{[\overline{\psi}_1,\overline{\psi}_2]}(\psi):= 
	\left\{
	\begin{array}{llllll}
		  1 & \ \ \text{for} \ \psi \in [\overline{\psi}_1,\overline{\psi}_2],\\
		  0 &  \ \ \text{otherwise}.
	\end{array}
	\right.
\end{equation*}
The computational times present in the numerical tests include the time of formation of the matrices, the setup time of the preconditioner (when used) and the solving time.

\subsection{Test with smooth solution} \label{sec:smooth}

We deal with a smooth solution on all the space--time domain $\Omega\times(0,1)$, with $\Omega=(0,1)$, to analyze the convergence order of SU method. 
We set $C_m=1$ and $\mathbf{D}(\vect{x})=10^{-4}\mathbb{I}_{N_s}$.

Moreover, in this test, we neglect the action of the recovery variable, i.e., we fix a null value for $w$ and we do not introduce the relaxations presented in Section \ref{sec:semiimplicit}.

We compute the function $f$ in \eqref{eq:McCulloch} in order to generate as exact solution $u_{\text{ex}}=10\sin(\pi x) \sin(\pi t) (1-\exp(-x))(1-\exp(x-1))(1-\exp(t-1))$. The resulting linear systems are solved by direct solver provided by {\sffamily MATLAB} (backslash operator ''\textbackslash'').

\begin{figure}[htbp]
	\centering
		\includegraphics[width=1.00\textwidth]{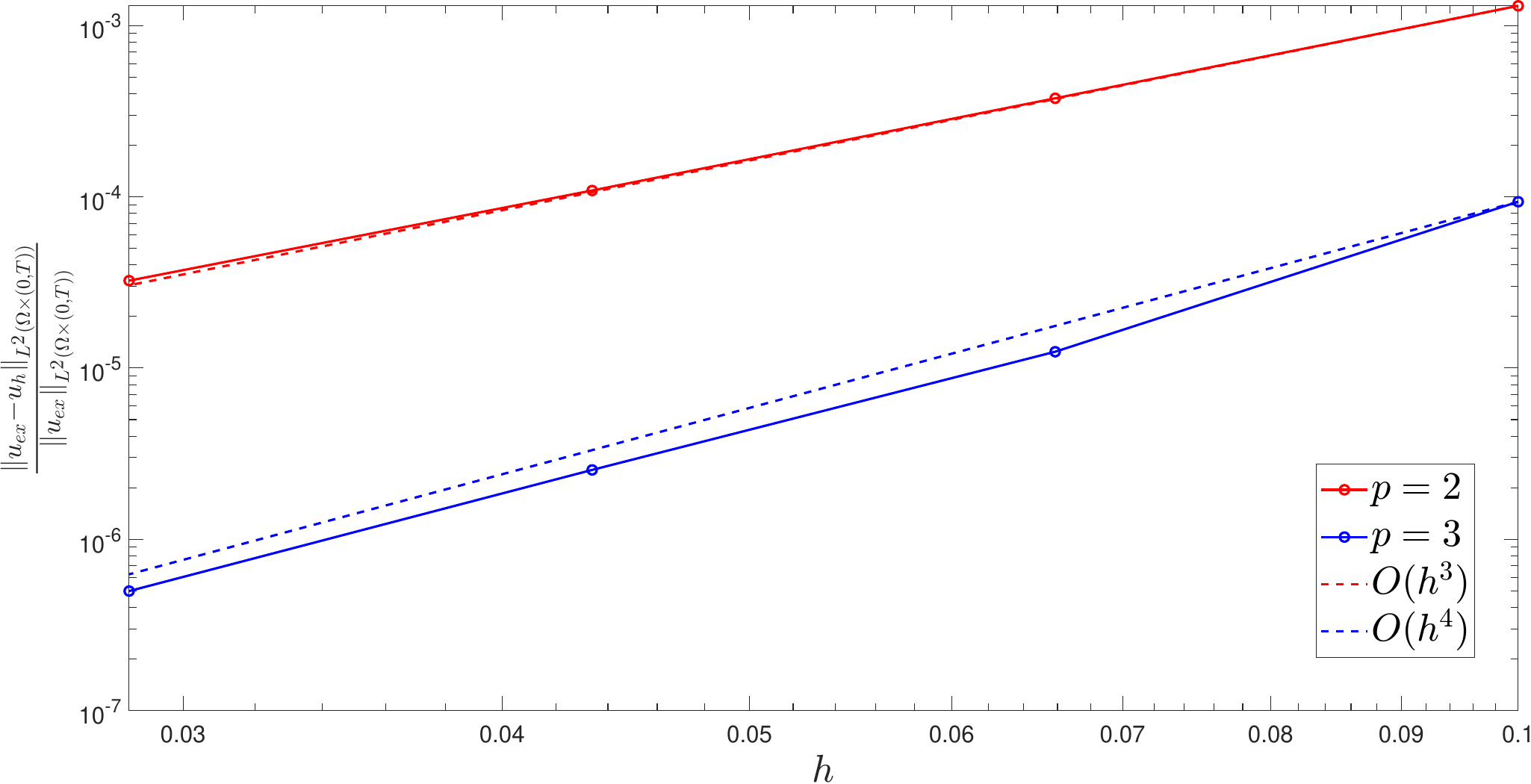}
	\caption{SU relative error plots in $L^2$-norm.}
	\label{fig:SU_ce_L2_error}
\end{figure}

Figure \ref{fig:SU_ce_L2_error} shows the error plot in \( L^2 \)-norm for SU method with uniform meshes in space and time, i.e., $h_1=h_t=:h$ and B-spline degrees  $p=2,3$, exhibiting optimal convergence order.

\subsection{2D spatial domain}\label{sec:McCulloch2D}
\begin{figure}[htbp]
	\centering
		\includegraphics[width=1\textwidth]{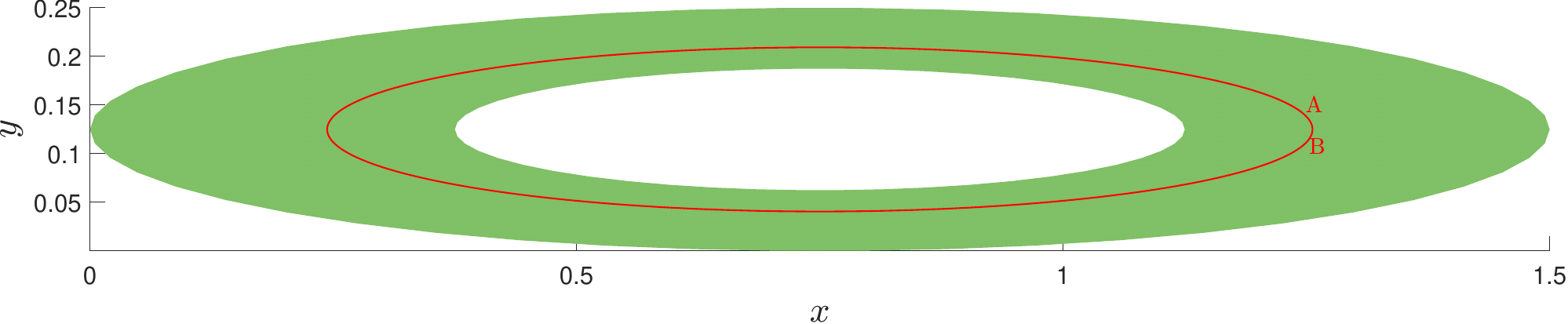}
	\caption{2D spatial domain with section line $s$ (A–B), test in Section \ref{sec:McCulloch2D}.}
	\label{fig:hole}
\end{figure}

In all of the tests present in this Section we consider as computational domain $\Omega$ an ellipse with semi-major axis equal to  $ 0.75$ and semi-minor axis equal to $ 0.125$, with an elliptic hole with semi-major axis equal to $0.375$ and semi-minor axis equal to $0.0625$. Figure \ref{fig:hole} shows an image of the spatial computational domain. The final time is  set equal to $T=300$. 
The local membrane capacitance and the anisotropy tensor are chosen as   $C_m=1$ and $\mathbf{D}(\vect{x})=10^{-4}\mathbb{I}_{N_s}$, respectively. 
The source term $f$ is as follows:
\begin{equation*} 
		f(x,y,t)=1/4 \, \exp\left(-\left(5\cdot 10^2\left(\left(y-L_2/2 \right)^2+\left(x-8/15 \, L_1/T\,t\right)^2\right)\right)\right) \ \chi_{[90,100]}(t),
\end{equation*}
where $L_2=0.125$ and $L_1=1.5$.

 The   linear systems are solved by direct solver provided by {\sffamily MATLAB} (backslash operator “\textbackslash'').

\subsubsection{SU stabilization with uniform and local refinement in time: comparison with Galerkin method}
\label{sec:2Dunif_local}
In this subsection we compare the results obtained with  SU stabilization with the ones obtained with classical Galerkin method, i.e., without the term $\widetilde{S}$ in \eqref{eq:solverSU_approx}.

\begin{figure}[htbp]
	\centering
		\centering
			\begin{subfigure}{1\textwidth}
				\includegraphics[width=1\textwidth]{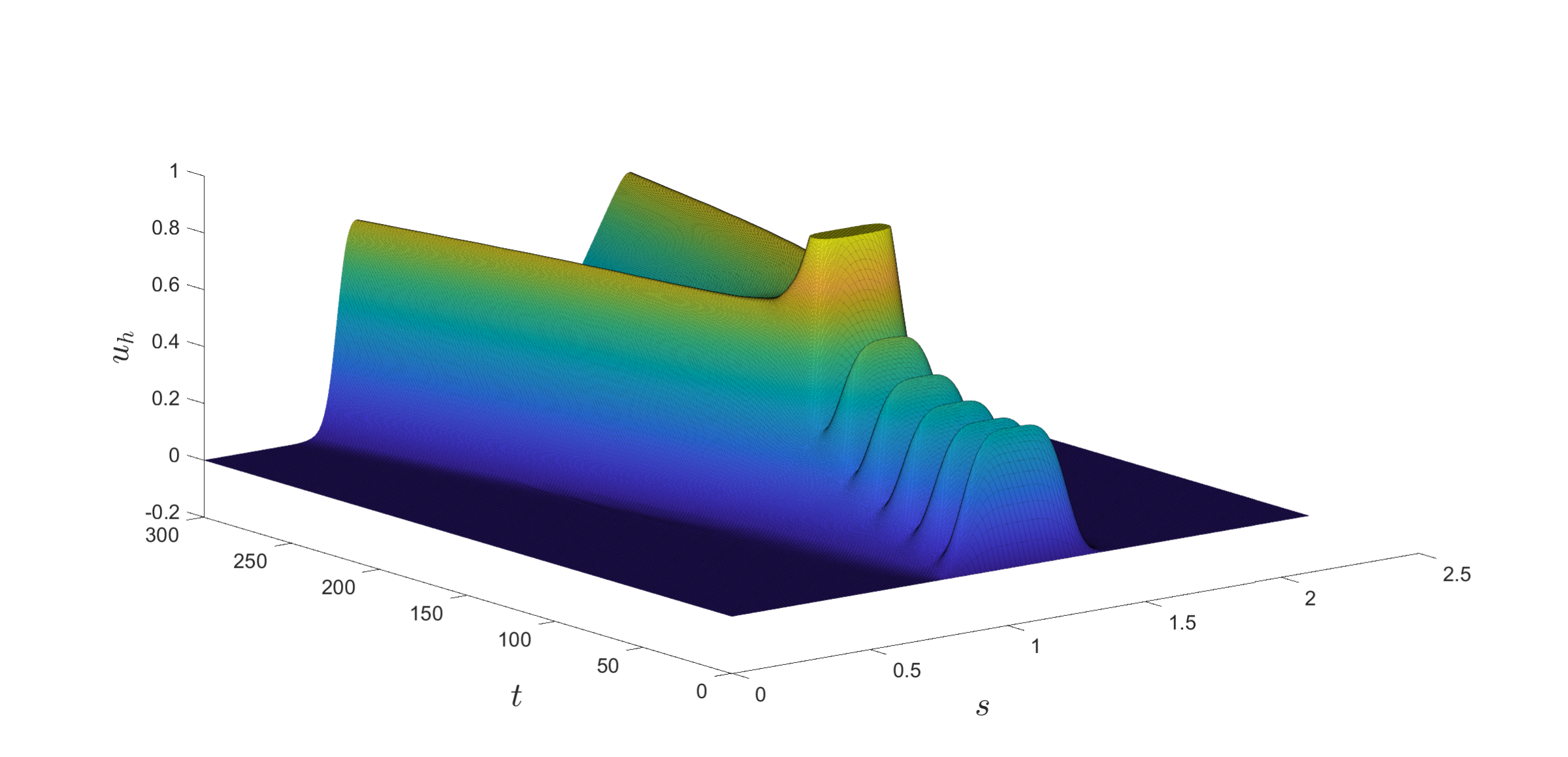}
			\end{subfigure}
			\begin{subfigure}{1\textwidth}
				\includegraphics[width=1\textwidth]{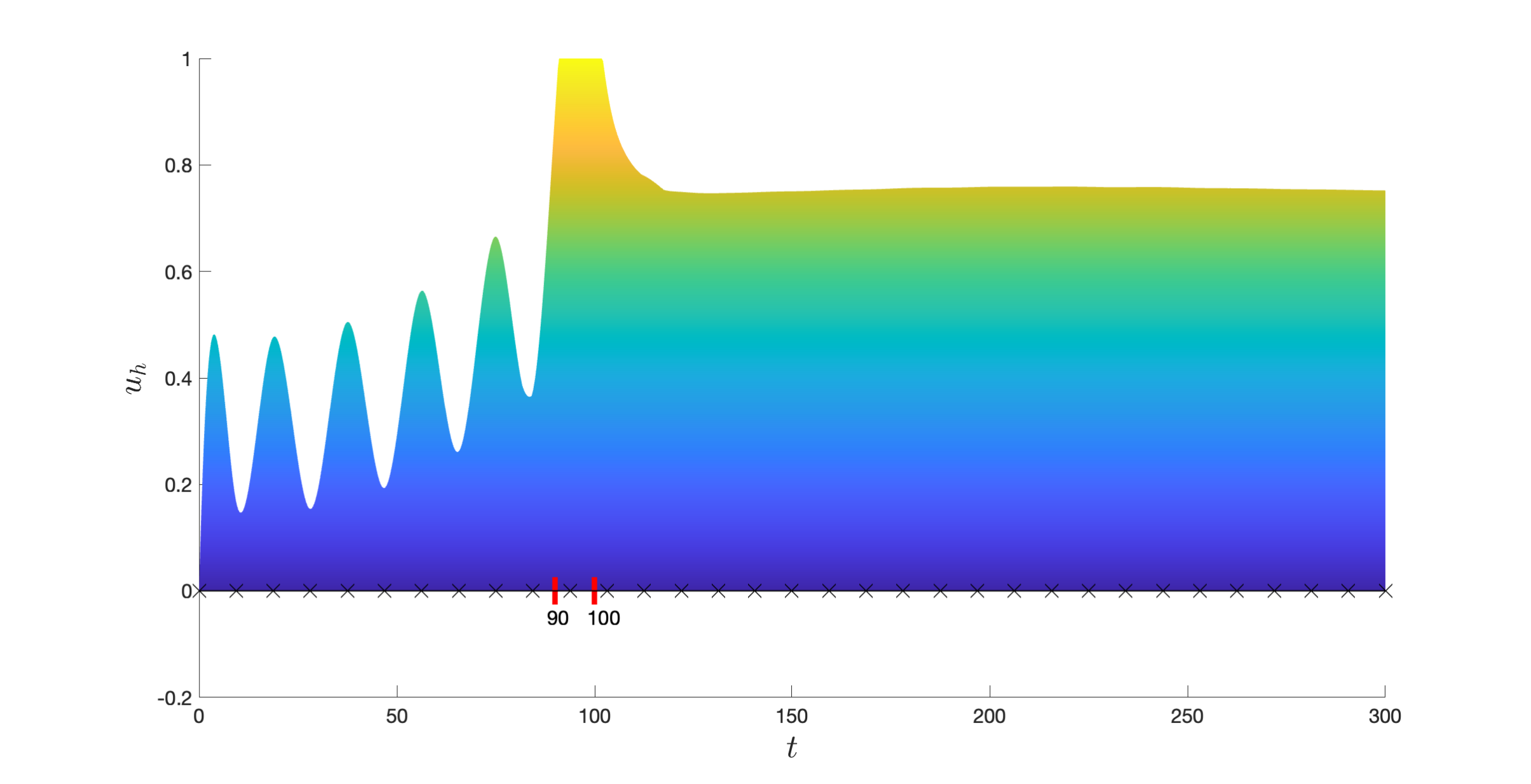}
			\end{subfigure}
	\caption{Two views of standard Galerkin transmembrane potential $u_h$ along section line $s$ on the 2D spatial domain for the  test in Section \ref{sec:McCulloch2D} (the extreme values of the source interval in time $[90,100]$ are highlighted in  red), with a uniform mesh in time represented with black crosses. The number of degrees of freedom  is $n_1=73, n_2=11$ and $N_t=34$.}
	\label{fig:McCulloch_hole_Galerkin}
\end{figure}

\begin{figure}[htbp]
	\centering
		\begin{subfigure}{1\textwidth}
			\includegraphics[width=1\textwidth]{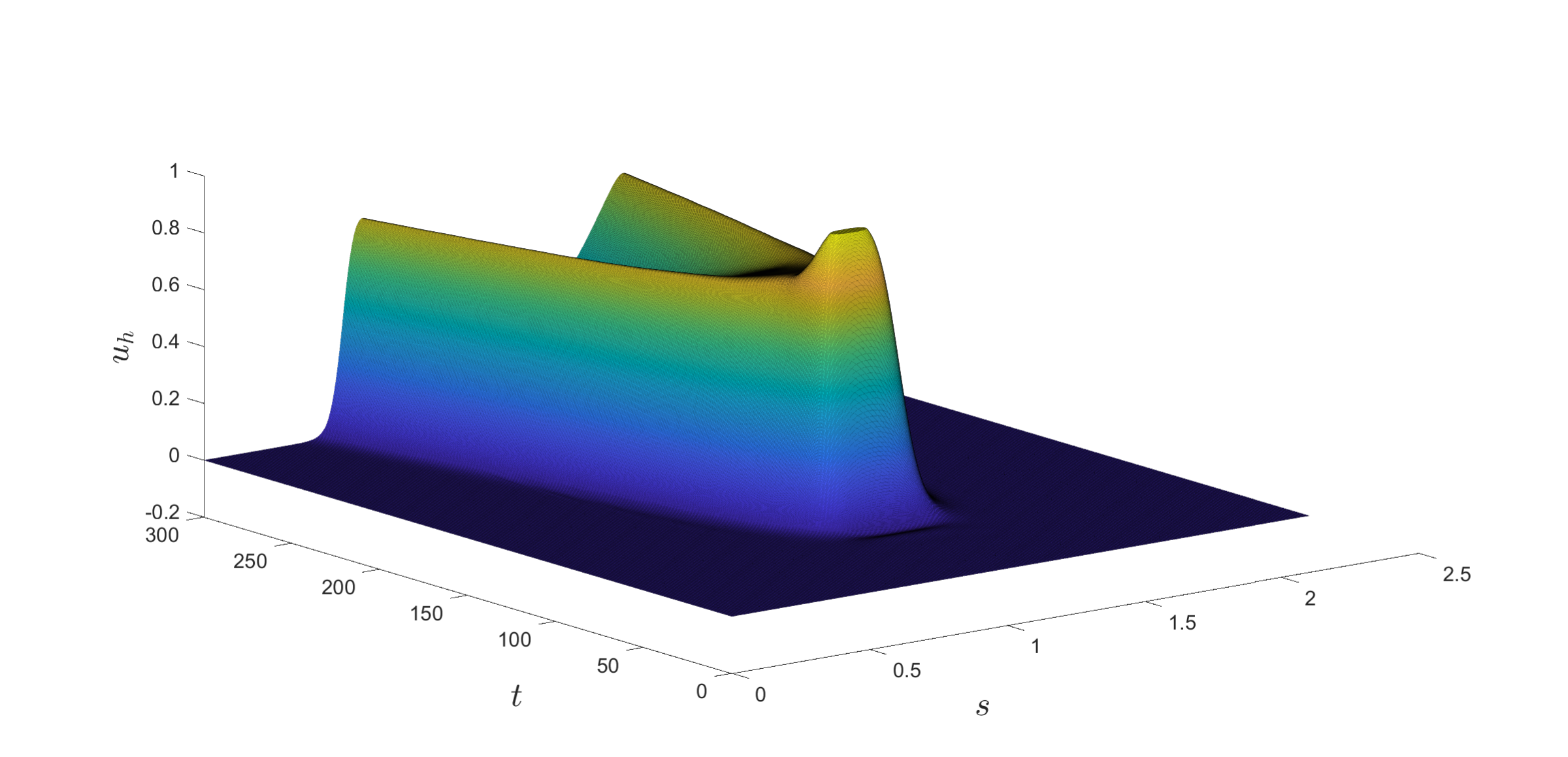}
		\end{subfigure}
		\begin{subfigure}{1\textwidth}
			\includegraphics[width=1\textwidth]{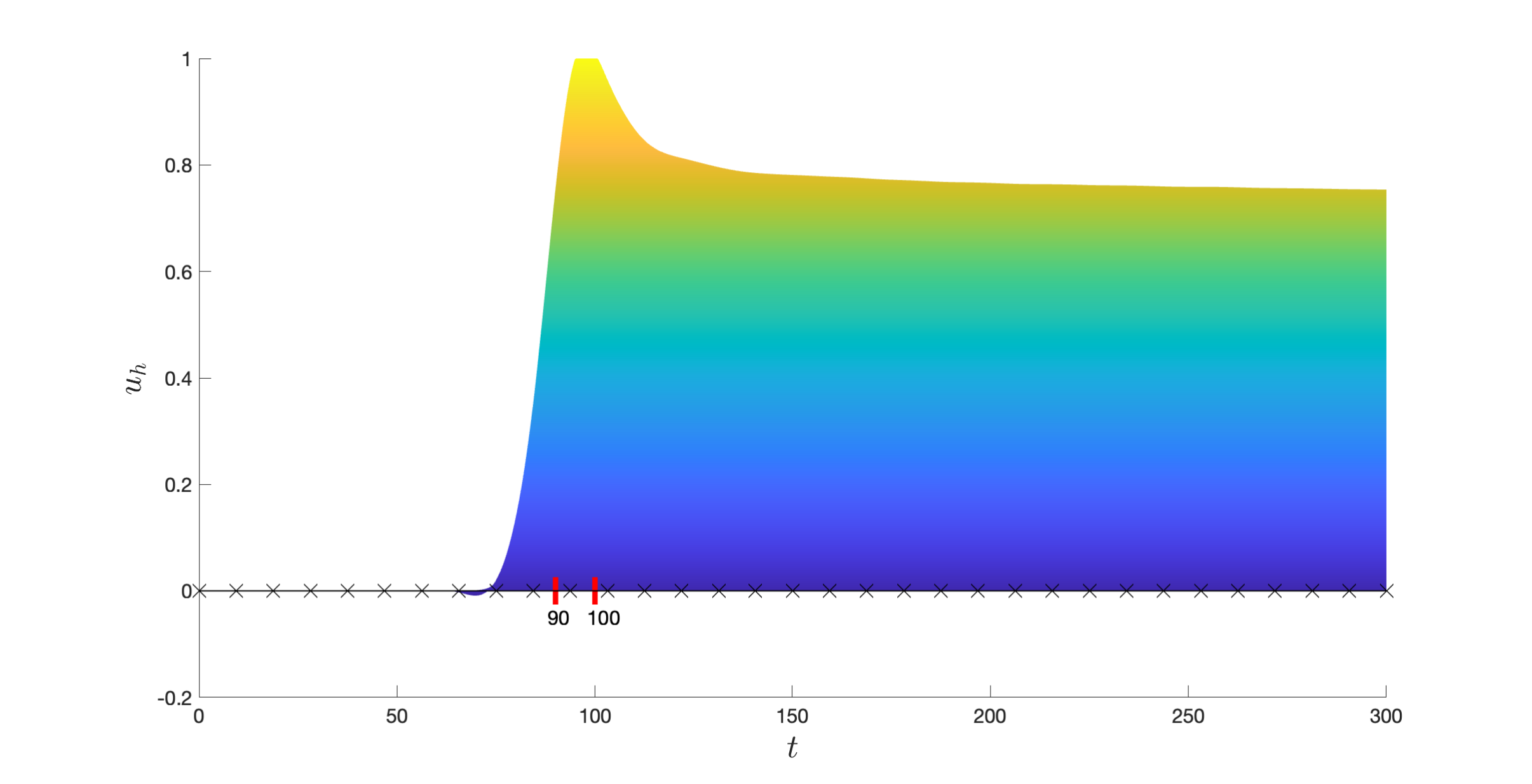}
		\end{subfigure}
	\caption{Two views of SU transmembrane potential $u_h$ along section line $s$ on the 2D spatial domain for the test in Section \ref{sec:McCulloch2D} (the extreme values of the source interval in time  $[90,100]$ are highlighted in red), with a uniform mesh in time represented with black crosses. The number of degrees of freedom  is   $n_1 = 73, n_2=11$ and $N_t=34$.}
	\label{fig:McCulloch_hole_SU}
\end{figure}

\begin{figure}[htbp]
	\centering
		\includegraphics[width=1\textwidth]{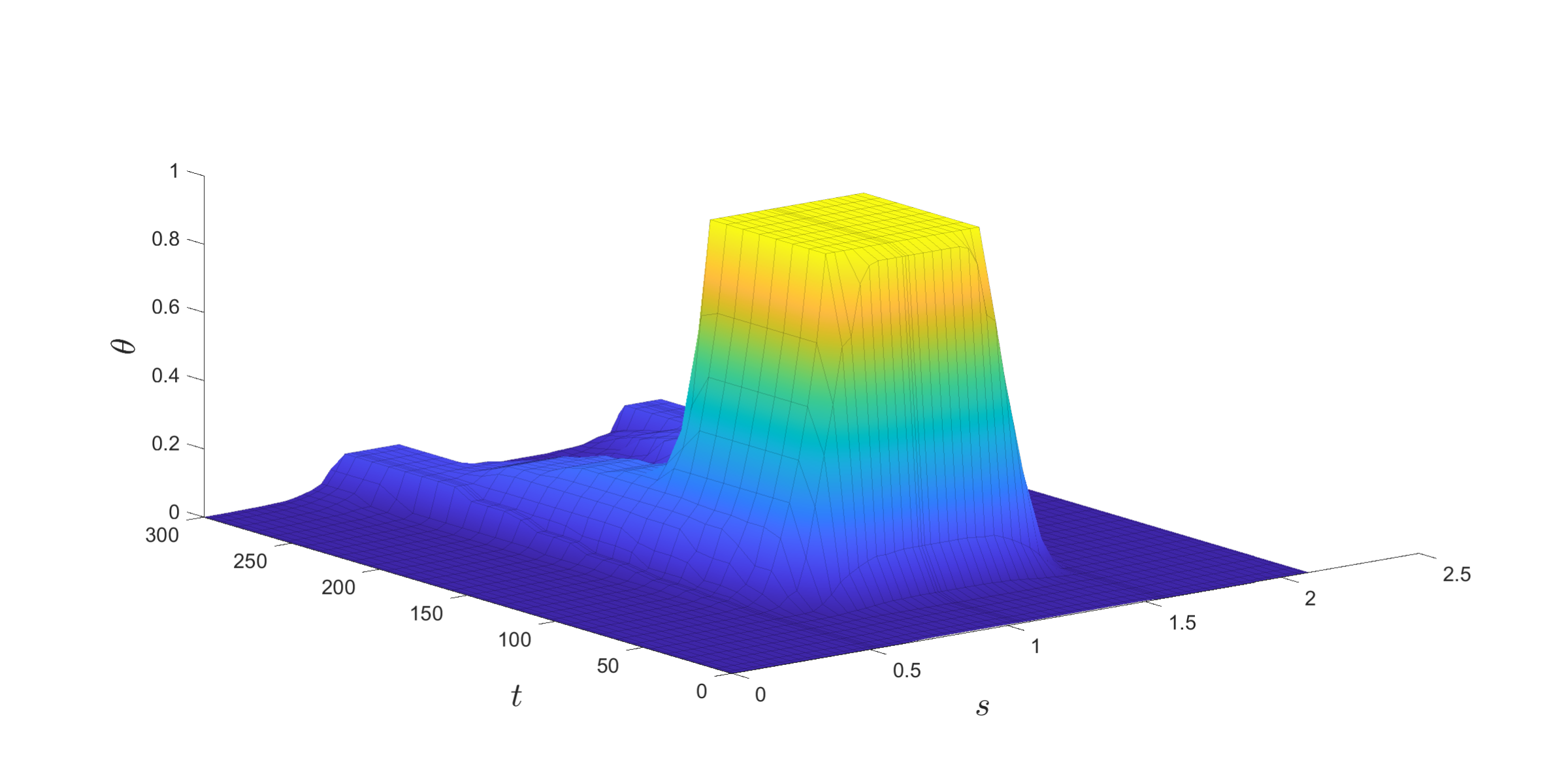}
	\caption{The function $\theta$ evaluated at the solution $(u_h,w_h)$ along section line $s$ on 2D spatial domain, test in Section \ref{sec:McCulloch2D}.}
	\label{fig:theta_HOLE}
\end{figure}

\begin{figure}[htbp]
	\centering
		\includegraphics[width=0.5\textwidth]{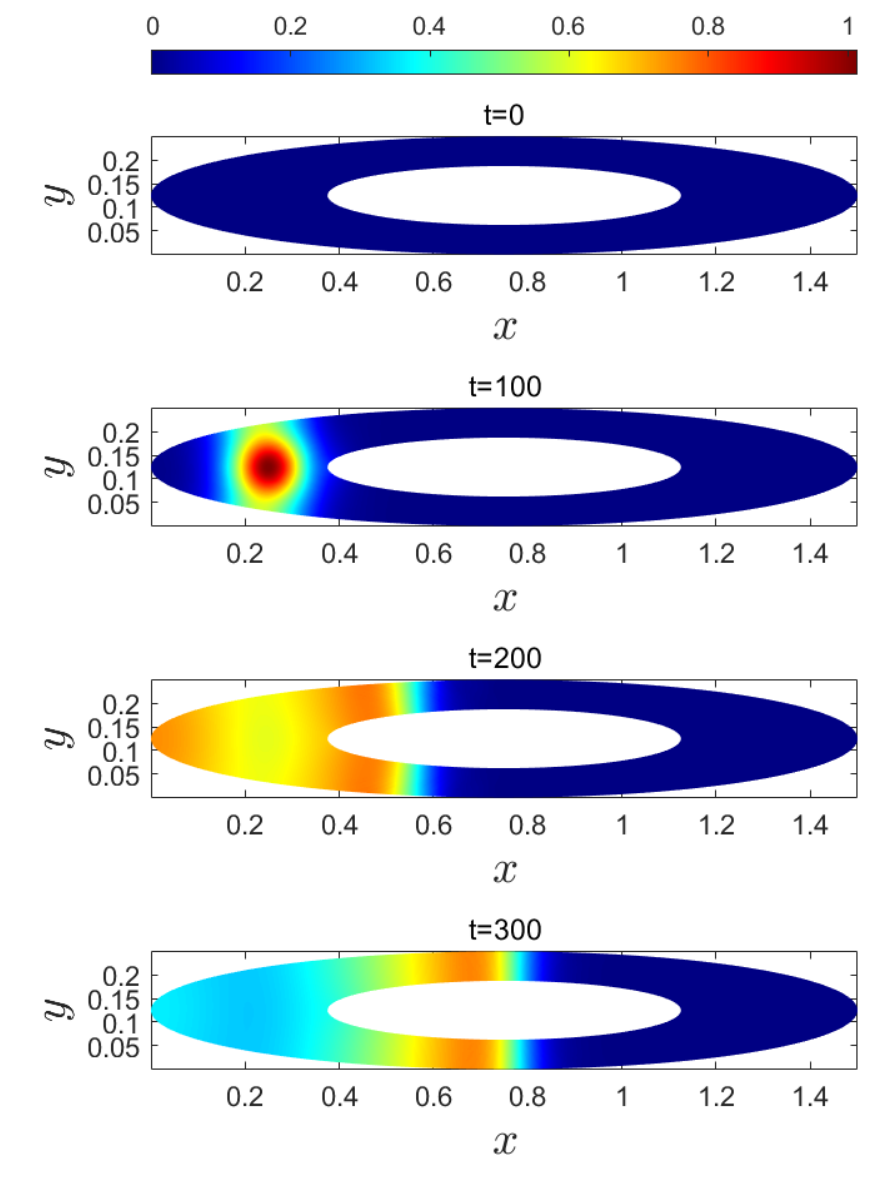}
	\caption{SU transmembrane potential, for various fixed times on the 2D spatial domain for the test in Section \ref{sec:McCulloch2D}.}
	\label{fig:McCulloch_hole_SU_T}
\end{figure}

For the first test we consider a uniform discretization with meshsizes $h_1=2^{-6}$, $h_2=2^{-3}$ and $h_t=2^{-5}$ and degree $p=3$.
Figures \ref{fig:McCulloch_hole_Galerkin} and \ref{fig:McCulloch_hole_SU} show the numerical solutions along section line $s$ (depicted in Figure \ref{fig:hole}) of the standard Galerkin and SU methods, while the function $\theta$ evaluated at the solution $(u_h,w_h)$, regulating the Spline Upwind activation, is presented in Figure \ref{fig:theta_HOLE}. Furthermore, Figure \ref{fig:McCulloch_hole_SU_T} shows the SU method solutions for various fixed times.
In this test the standard Galerkin method presents spurious oscillation and these numerical instabilities are reduced using the SU method. We remark that with SU method the spurious behaviour before the activation time  $t=90$  is contained in the support extension of a single basis function and thus reducing the support extent, e.g., by refining near the activation source, could significantly improve the performance of the method.

Table \ref{tab:timehole} shows that the computational time of the standard Galerkin method is considerably higher compared to that of the SU method. This result is due to the high number of iterations required to reach the set convergence, caused by considerable numerical instabilities.

\begin{table}[htbp]
	\centering
		\begin{tabular}{|l|c|r|}
			\hline
			 & Fixed point iterations & Time (s)  \\
			\hline
			Standard Galerkin method & 145 & $1.9\cdot 10^3$ \\
			\hline
			SU method & 31 & $4.6 \cdot 10^2$\\
			\hline
		\end{tabular}
	\caption{Computational cost comparison for 2D spatial domain with uniform knot vector, test in Section \ref{sec:McCulloch2D}.}
	\label{tab:timehole}
\end{table}

\begin{figure}[htbp]
	\centering
			\includegraphics[width=1\textwidth]{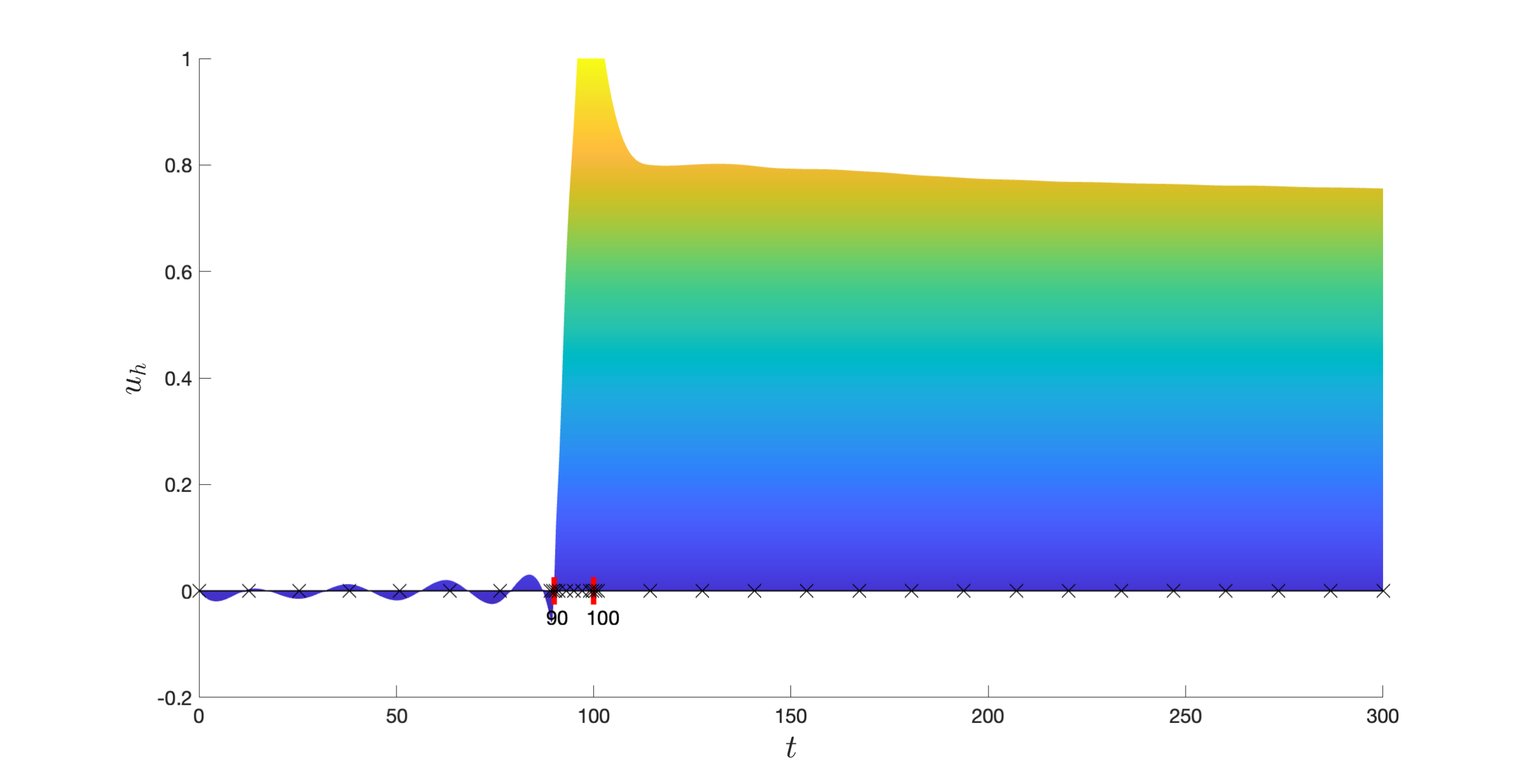}
	\caption{A view of standard  Galerkin transmembrane potential $u_h$ along section line $s$ on 2D spatial domain, test in Section \ref{sec:McCulloch2D} with a non-uniform mesh in time represented with black crosses. The number of degrees of freedom   is $n_1= 73,  n_2=  11$ and $  N_t= 34$.}
	\label{fig:McCulloch_hole_Gal_nonunif}
\end{figure}

\begin{figure}[htbp]
	\centering
			\includegraphics[width=1\textwidth]{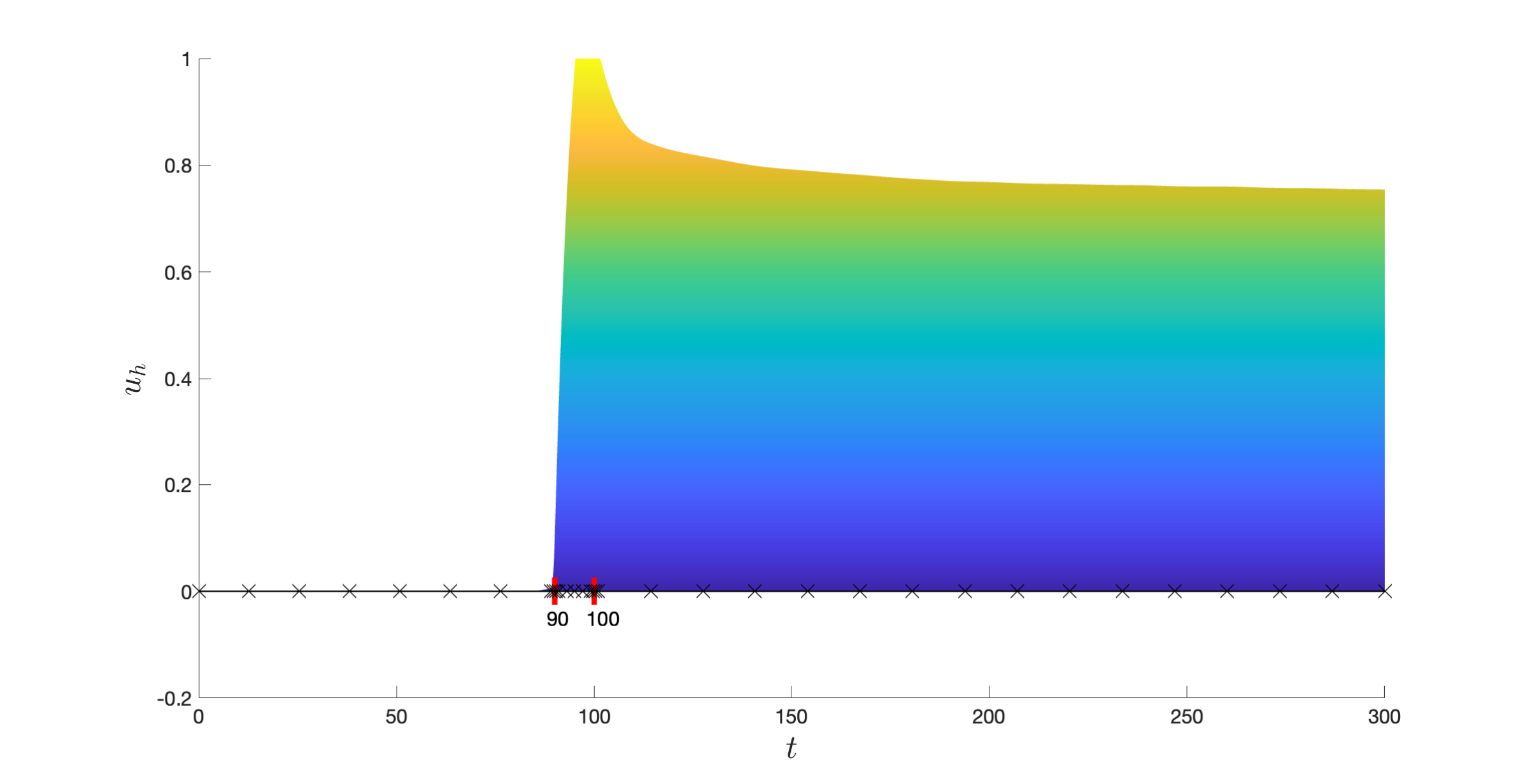}
	\caption{A view of SU transmembrane potential $u_h$ along section line $s$ on 2D spatial domain, test in Section \ref{sec:McCulloch2D} with a non-uniform mesh in time represented with black crosses. The number of degrees of freedom    is  $n_1= 73,  n_2=  11$ and $  N_t= 34$.}
	\label{fig:McCulloch_hole_SU_nonunif}
\end{figure}

We remark that we set up the test case in an unfavorable condition from the perspective of stability by employing a uniform knot partition. We can further improve the behavior of SU stabilization exploiting the knowledge of the location of the source function by refining in the time direction in the regions  of activation and shutdown  of the source $f$, to better capture the discontinuities. In particular, we perform a second test   keeping fixed the degrees and  the number of degrees of freedom and using a non-uniform knot vector in time direction, represented in Figure \ref{fig:McCulloch_hole_Gal_nonunif}  and \ref{fig:McCulloch_hole_SU_nonunif} with black crosses, with  $p+1$ knots clustered around $t=90$ and $p+1$ knots clustered around $t=100$ and a almost uniform partition elsewhere. While the oscillations are still present with Galerkin method, the SU method greatly benefits from the locally refined in time partition, as shown in Figure  \ref{fig:McCulloch_hole_Gal_nonunif} and Figure  \ref{fig:McCulloch_hole_SU_nonunif}, respectively. The support extent of the basis is reduced near the activation source, resulting in improved results. We also remark that the use of local mesh refinement techniques,   as \cite{Feischl2016,Vuong2011}, could automatize the refinement procedure in the source region.

\subsubsection{Comparison with low regularity $C^0$ splines in time}
In this section we present   numerical results obtained with low regularity ($C^0$) splines in the time direction. As our stabilization technique can not be extended to splines of regularity lower than $C^{p-1}$, we use a classical stabilization technique: namely, the SUPG method \cite{Brooks1982}. Moreover, to enhance the stability of SUPG, following  \cite{Bazilevs2007}, we add a shock capturing term. In particular, we replace $ \mathcal{S}$ in \eqref{eq:McCulloch_SU} with $  {\mathcal{S}}_{SUPG} + {\mathcal{S}}_{SC}$
  where
  $$ 
\mathcal{S}_{SUPG}({u}_h,v_h,w_h):=\frac{h}{2p+1}\int_0^T\int_{\Omega}     \left(C_m\partial_t {u_h}-   \nabla\cdot(\mathbf{D}\nabla  {u}_h)  +C_r(u_h, w_h) {u_h}-f\right)\partial_t v_h\ \d\Omega\dt
$$
and
 $$ 
\mathcal{S}_{SC}({u}_h,v_h,w_h):=c\frac{h^2}{4}\int_0^T\int_{\Omega}  \frac{|R(u_h,w_h)|}{C_m(T^{-1}\|u_h\|_{\infty}+\|\partial_t u_h\|_{\infty})} \partial_t {u}_h\partial_t v_h\ \d\Omega\dt,
$$
and where $R(u_h,w_h):=C_m\partial_tu_h-\nabla\cdot(\mathbf{D}\nabla u_h)+c_1u_h(u_h-a)(u_h-1)+c_2u_hw_h-f$ is the residual.
The value $c\geq 0$ is a tuning parameter.
 Then the fixed point scheme \eqref{eq:solverSU} is consequently modified by replacing 
 ${\mathcal{S}}(\tilde{u}_h^{k+1},v_h)$ with $ {\mathcal{S}}_{SUPG}(\tilde{{u}}^{k+1},u_h^k,v_h,w_h^k)+ {\mathcal{S}}_{SC}(\tilde{{u}}^{k+1},u_h^k,v_h,w_h^k)$, where  with a little abuse of notation, we defined
  \begin{align*} 
\mathcal{S}_{SUPG}(\tilde{{u}}^{k+1},u_h^k,v_h,w_h^k):=\frac{h}{2p+1}&\int_0^T\int_{\Omega}    (  C_m\partial_t {\tilde{{u}}^{k+1}}-   \nabla\cdot(\mathbf{D}\nabla \tilde{{u}}^{k+1})  +C_r( {{u}}^{k}, w_h^k) \tilde{{u}}^{k+1}-f)\partial_t v_h\ \d\Omega\dt
\end{align*}
and
 \begin{align*} 
\mathcal{S}_{SC}(\tilde{{u}}^{k+1},u_h^k,v_h,w_h^k):=c\frac{h^2}{4}\int_0^T\int_{\Omega}  \frac{|R( {u}^h_h,w^k_h)|}{C_m(T^{-1}\| {u}^k_h\|_{\infty}+\|\partial_t  {u}^k_h\|_{\infty})} \partial_t \tilde{u}^{k+1}_h\partial_t v_h\ \d\Omega\dt.
\end{align*}
 Then, to speed up computations and following the same steps to define $\widetilde{S}$ in Section 3.1, we  define  $\widetilde{\mathcal{S}}_{SUPG}(\tilde{{u}}^{k+1},u_h^k,_h,v_h,w_h^k)$ and $\widetilde{\mathcal{S}}_{SC}(\tilde{{u}}^{k+1},u_h^k,_h,v_h,w_h^k)$ by taking ${\mathcal{S}}_{SUPG}(\tilde{{u}}^{k+1},u_h^k,_h,v_h,w_h^k)$ and $ {\mathcal{S}}_{SC}(\tilde{{u}}^{k+1},u_h^k,_h,v_h,w_h^k)$, respectively, and  replacing $C_r(u_h^k,w_h^k)$ and $ \frac{|R( {u}^h_h,w^k_h)|}{C_m(T^{-1}\| {u}^k_h\|_{\infty}+\|\partial_t  {u}^k_h\|_{\infty})} $  with  low-rank approximations. The fixed point scheme to be solved is then \eqref{eq:solverSU_approx} with $ \widetilde{\mathcal{S}}_{SUPG}+\widetilde{\mathcal{S}}_{SC}$ that replaces  $\widetilde{\mathcal{S}}.$ 
 
\begin{center}
\begin{figure}[H]
\includegraphics[scale=0.4]{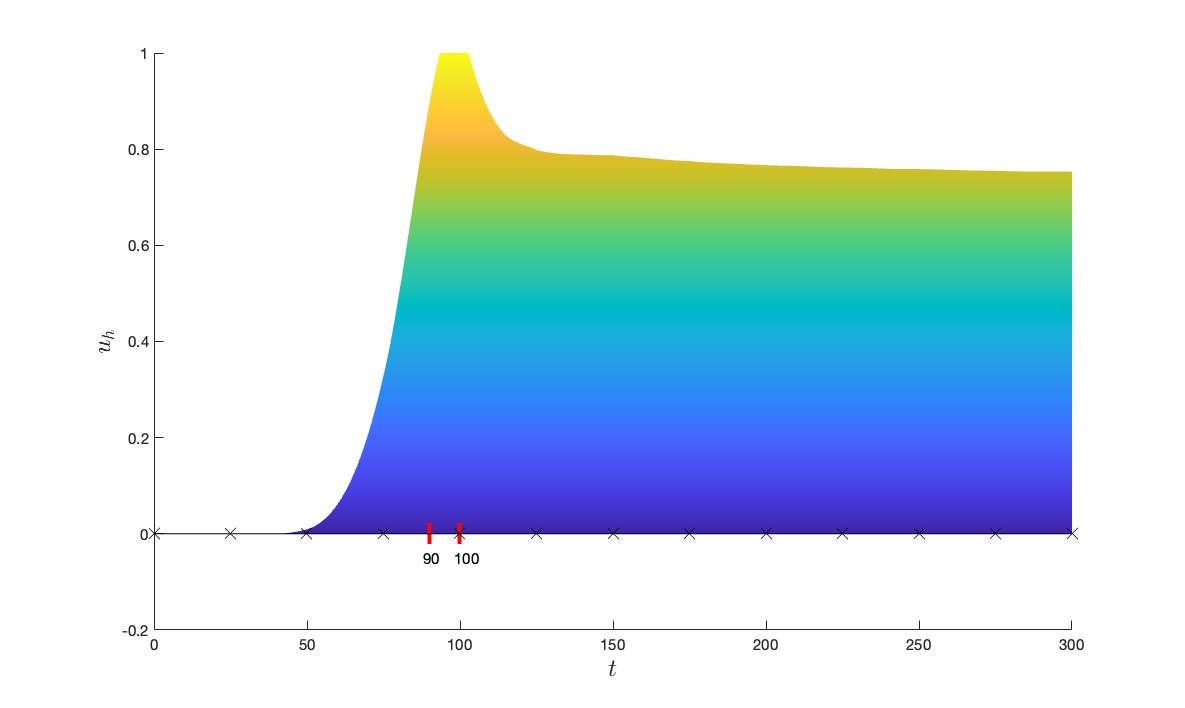}
\caption{A view of SU transmembrane potential $u_h$ along section line $s$ on 2D spatial domain, test in Section \ref{sec:McCulloch2D} with a vuniform mesh in time represented with black crosses. $C^0$ continuous splines in time and SUPG stabilization with shock capturing are used. The number of degrees of freedom is $n_1=73$, $n_2=11$ and $N_t=36$.}
\label{fig:C0_supg_sc}
\end{figure}
\end{center}

In our test, we consider $p=3$, $C^{p-1}$ continuity in all the spatial directions and $C^0$ continuity in the time direction. We take  the same meshsizes in space as the ones used in Section \ref{sec:2Dunif_local} to generate Figure \ref{fig:McCulloch_hole_SU}, i.e.,  $h_1=2^{-6}$ and $h_2=2^{-3}$,  and a  meshsize in the time direction that makes the number of degrees of freedom in time comparable to the ones used for the SU test with uniform refinement in Section \ref{sec:2Dunif_local}. In particular, we choose $h_t=1/12$, resulting in  $N_t=36$ degrees of freedom in time. We remark  that this mesh is coarser than the one used to generate Figure \ref{fig:McCulloch_hole_SU}. 
 The results are shown in Figure \ref{fig:C0_supg_sc}. Note that the spurious behavior before the activation time  $t=90$  is contained in the support extension of a single basis function, exactly as it happens as with the SU stabilization (see Figure \ref{fig:McCulloch_hole_SU}). The SUPG stabilization with shock capturing employs 106 iterations and $1.87\cdot 10^3$ seconds to converge (cfr. Table \ref{tab:timehole}).
Comparing Figure \ref{fig:C0_supg_sc} with Figure \ref{fig:McCulloch_hole_SU}, we conclude that, for a fixed number of degrees of freedom, SU stabilization  seems to be more effective than SUPG with shock capturing: the area of spurious oscillations  is smaller when the former method is used.


\subsection{3D spatial domain}\label{sec:Ventr}
\begin{figure}[htbp]
	\centering
	\includegraphics[width=0.7\textwidth]{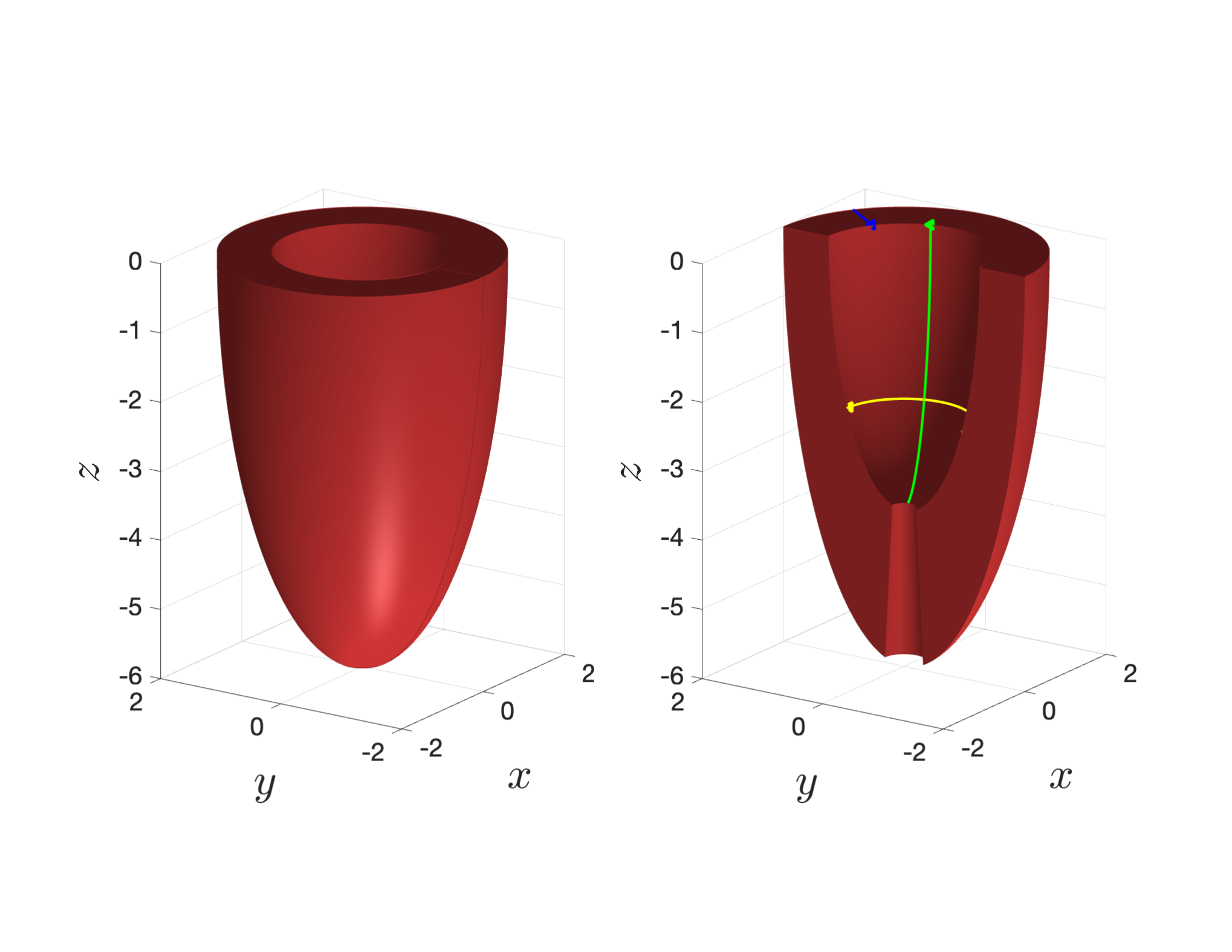}
	\caption{Ventricular spatial domain $\Omega$ for the test in Section \ref{sec:Ventr} on the left, and its corresponding cross-section with arrows representing the mapping of parametric coordinates ($\eta_1$ in yellow, $\eta_2$ in green, $\eta_3$ in blue).}
	\label{fig:3D_Ventr_domain}
\end{figure}

In this test case, we consider as 3D-spatial domain a commonly used idealization of the left ventricle: a truncated ellipsoid, as represented in Figure~\ref{fig:3D_Ventr_domain}. For the parametrization of the ellipsoid, we refer to \cite{Charawi2017}. 
This is a significant benchmark, often employed to assess the good behavior of the  solving methods for this kind of problems (see   e.g. \cite{Willems2024,Charawi2017,franzone1998spread} and references therein).
We fix $T=300$ and $p=2$. 

We set $C_m=1$ and   we consider $h_1=h_2=h_3=2^{-5}$ and $h_t=2^{-6}$.   We suppose to have axial isotropy, i.e., same conductivity in both tangential and normal directions, and that the fibers rotate intramurally linearly with the depth for a total amount of 120 degrees proceeding counterclockwise from epicardium and endocardium. We thus consider the conductivity tensor $\mathbf{D}(\vect{x})$ defined as 
$$
\mathbf{D}(\vect{x})=\frac{\lambda}{1+\lambda}\left(\sigma_t\mathbb{I}_{N_s}+(\sigma_l-\sigma_t)\mathbf{a}(\vect{x})\mathbf{a}(\vect{x})^T\right)
$$ 
where $\lambda= {\sigma_l}\slash{\sigma^e_l}$ and $\sigma_l=3\times 10^{-3}, \sigma_t=3.1525\times 10^{-4}$ and $\sigma_l^e=2\times 10^{-3}$ are the conductivity coefficients and $\mathbf{a}$ is the unit vector tangent to the cardiac fiber at a point $\vect{x}$. 
For more details about the choice of the conductivity tensor we refer to \cite{pennacchio2011fast}.
 
The linear systems \eqref{eq:lin_sys_A} are solved using the generalized minimal residual method (GMRES) with a modified version of the preconditioner presented in Section \ref{sec:prec}, that incorporates some information on the geometry and the coefficients (see \cite[Section 4.4]{Loli2020} for details). To solve \eqref{eq:lin_sys_L}, we exploit the technique presented in Section \ref{sec:prec} and we solve \(N_s\) independent systems associated with the \(N_t \times N_t\) matrix with direct solver provided by {\sffamily MATLAB} (backslash operator ''\textbackslash''), while \(N_t\) independent systems associated with the \(N_s \times N_s\) mass matrix are solved by using the preconditioned conjugate gradient method (PCG)  with the  preconditioner described in Section \ref{sec:prec}.
We fix the tolerance for both methods equal to $10^{-8}$.

The source term $f$ is defined as $f(x,y,z,t) =  \hat{f}\circ\vect{G}^{-1}(x,y,z,t) $, where $\tilde{f}:\widehat{\Omega}\times [0,1]\rightarrow \mathbb{R}$ is 
\begin{equation*} 
	\hat{f}(\eta_1,\eta_2,\eta_3,\tau)=  \left\lbrace \begin{array}{ll} \chi_{[0.15,0.2]}(\tau) \quad & \text{if}\  (\eta_1-0.5)^2+(\eta_2-0.5)^2+\eta_3^2\leq 0.0025,\\
	 0 \quad & \text{otherwise},\end{array}\right.
\end{equation*} 
and $\vect{G}$ is the geometrical mapping (see Section 1). 

\begin{figure}[htbp]
	\centering
	\includegraphics[width=0.9\textwidth]{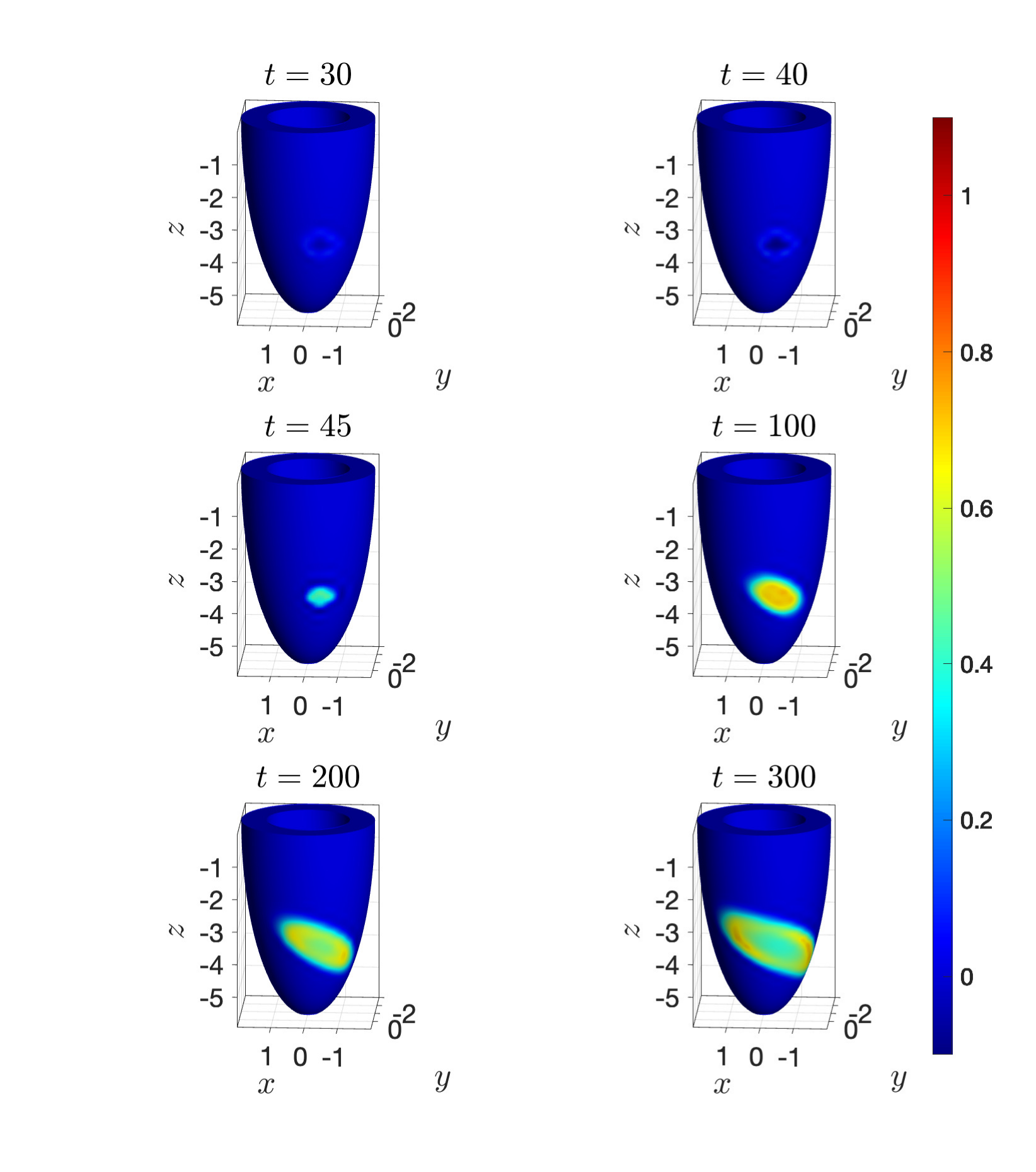}
	\caption{ Standard Galerkin transmembrane potential, for various fixed times on 3D spatial domain for test in Section \ref{sec:Ventr}.  The number of degrees of freedom  is $n_1= 37,  n_2=  34,   n_3= 34$ and $  N_t= 65$.}
	\label{fig:3D_Ventr_Galerkin}
\end{figure}

\begin{figure}[htbp]
	\centering
	\includegraphics[width=0.9\textwidth]{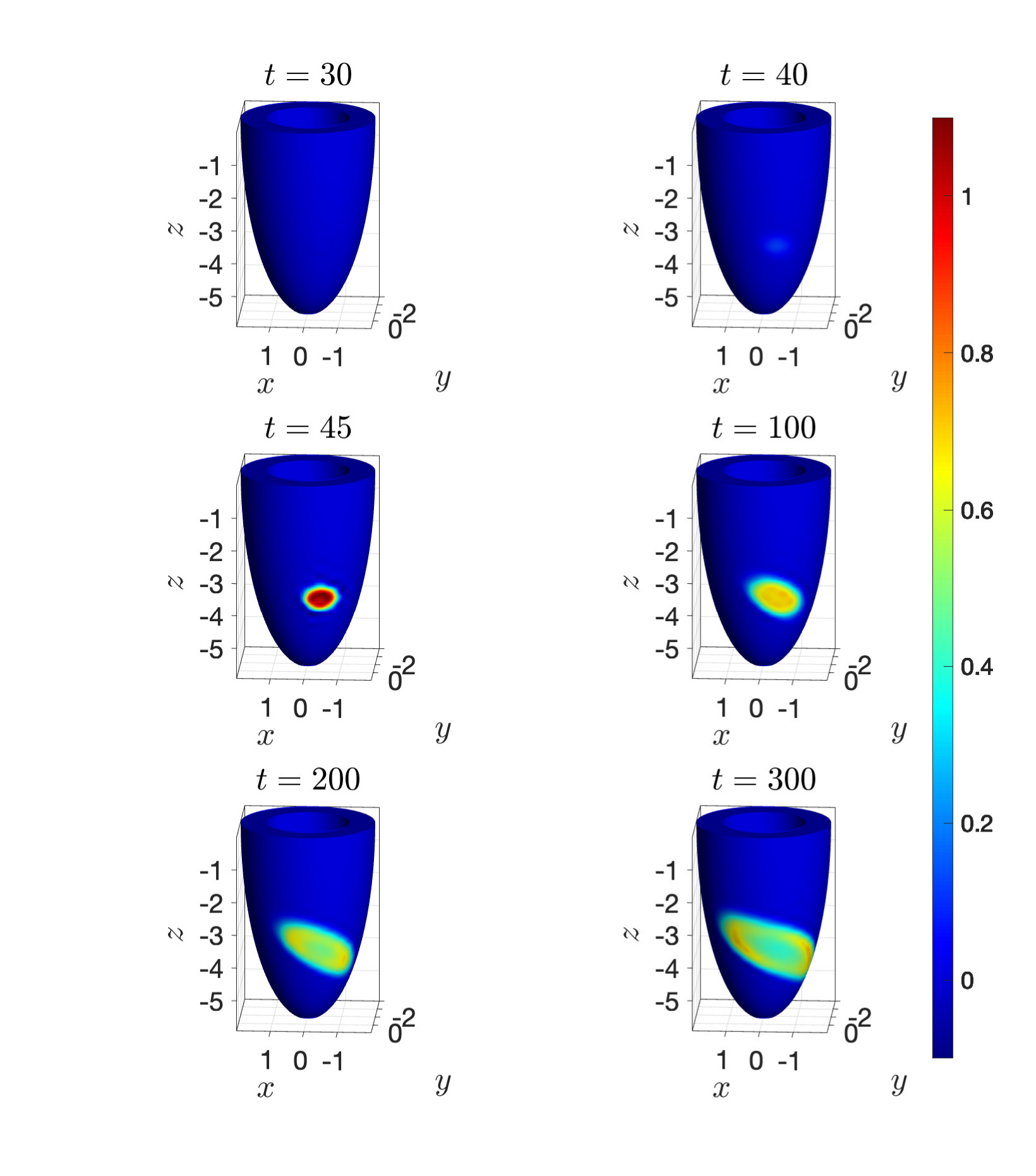}
	\caption{SU transmembrane potential, for various fixed times onthe  3D spatial domain for the  test in Section \ref{sec:Ventr}.  The number of degrees of freedom is $n_1= 37,  n_2=  34,   n_3= 34$ and $  N_t= 65$.}
	\label{fig:3D_Ventr_SU}
\end{figure}

Figures \ref{fig:3D_Ventr_Galerkin} and \ref{fig:3D_Ventr_SU} compare the numerical solutions obtained with the standard  Galerkin method and the SU method at various fixed times. Before the activation of the source term, i.e., for $t < 45$, where zero transmembrane potential is expected,   the standard Galerkin solution represented in Figure \ref{fig:3D_Ventr_Galerkin}  shows non-physical behaviors and at the activation time, i.e., for $t=45$, the solution does not reach the value of 1, indicating numerical instabilities. On the other hand, Figure~\ref{fig:3D_Ventr_SU} suggests that the SU method effectively reduces these oscillations: before the activation of the source term, the solution is almost zero in the spatial domain and it activates correctly at $t=45$.

Table \ref{tab:time3DVentr} presents the total computational time to solve the problem with the proposed method and the average  of GMRES and PCG   iterations for each fixed-point iteration. In particular, the SU method requires one and a half the total computation time of the standard Galerkin method.

\begin{table}[htbp]
	\centering
		\begin{tabular}{|l|c|c|r|}
			\hline
			 & GMRES / PCG average iterations & Fixed-point iterations & Time (s) \\
			\hline
			Standard Galerkin method & 143 / 4  & 43  & $ 8.2  \cdot 10^4$ \\
			\hline
			SU method &  149 / 4 & 44  & $1.3 \cdot 10^5$\\
			\hline
		\end{tabular}
	\caption{Computational cost comparison for 3D spatial domain, test in Section \ref{sec:Ventr}.}
	\label{tab:time3DVentr}
\end{table}


\section{Conclusions}\label{sec:conclusion}

In this work,  we proposed a solver to prevent spurious oscillations for evolutionary problems related to cardiac electrophysiology in the space--time IgA framework. We used a variant of the Spline Upwind (SU) method introduced in \cite{Loli2023}.

The numerical tests highlight the optimal convergence order for smooth solution and the stable behavior in the presence of sharp layers in the monodomain equation with the Rogers--McCulloch ionic model, both in two-dimensional and three-dimensional spatial domains. 

To complete the analysis, we quantified the computational cost in terms of the number of iterations required by the nonlinear solvers to obtain the solution, as well as the corresponding time. In our tests, transitioning from the standard Galerkin method to the SU method, we have showed that the computational times of the SU method are competitive with those of the standard Galerkin method. This makes SU stabilization appealing not only in terms of solution stability but also in terms of computational efficiency.  Moreover, refinement in time direction localized near the activation place of the source function improves the quality of the proposed stabilization method. 

While modeling the left ventricle as a fixed geometry is only  an approximation, it captures the essential features of cardiac electrical activity.
 Coupling more physics yields a much more complex problem that is definitely interesting to address thorough the proposed IgA technique. The first works in that direction are  \cite{torre2022efficient,willems2023isogeometric}. The challenges in extending our space-time approach to tackle these problems are many: the tensoriality between space and time is lost, the extension of the proposed stablization is not trivial and the computational cost of assembling and solving the resulting linear systems needs to be carefully investigated.
 
As  further extensions, we will consider also the bidomain model, which provides a more detailed representation of the electrical behavior of cardiac tissue by explicitly modeling the interaction between the intracellular and extracellular spaces. This extension is not trivial, as  the equations in this model are strongly coupled and it depends on the formulation that is used. This study deserves a deeper investigation.


\section*{Acknowledgements}
The authors wish to thank Dr.~A.~Bressan and  Prof.~S.~Scacchi for helpful preliminary discussions.\\
P.~Tesini is supported by European Union - Next Generation EU via the PRIN project ``ASTICE''.\\
L.~Dedè acknowledges the support of the PRIN project 202232A8AN ``Computational modeling of the human heart: from efficient numerical solvers to cardiac digital twins", MUR, Italy.\\
G.~Loli and G.~Sangalli acknowledge support from PNRR-M4C2-I1.4-NC-HPC-Spoke6 and the PRIN 2022 PNRR project NOTES (No. P2022NC97R).\\
M.~Montardini acknowledges support from the PRINN 2022 PNRR project HEXAGON (No. P20227CTY3).\\
P. F. Antonietti acknowledges partial support from ICSC - Centro Nazionale di Ricerca in High Performance Computing, Big Data, and Quantum Computing funded by European Union - NextGenerationEU.\\
The authors are members of the Gruppo Nazionale Calcolo Scientifico - Istituto Nazionale di Alta Matematica (GNCS-INdAM).\\
G.~Loli, M.~Montardini and P.~Tesini are partially supported by INdAM-GNCS Project ``Sviluppo di metodi numerici innovativi ed efficienti per la risoluzione di PDE''.
\begin{figure}[H]
	{\centering
		\hfill\includegraphics[scale=0.075]{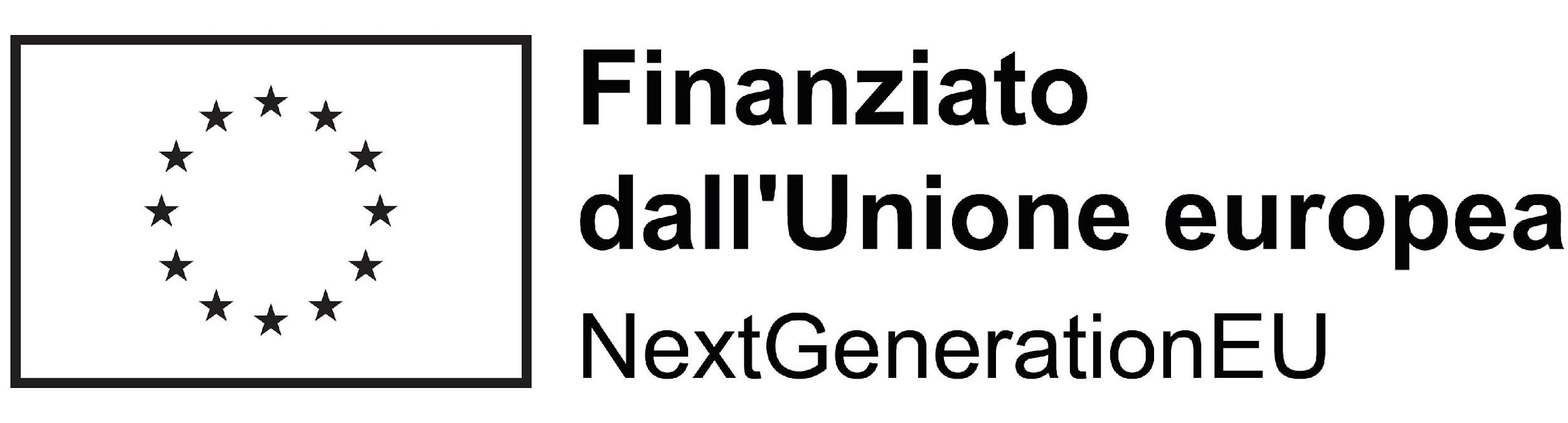}
		\hfill\includegraphics[scale=0.075]{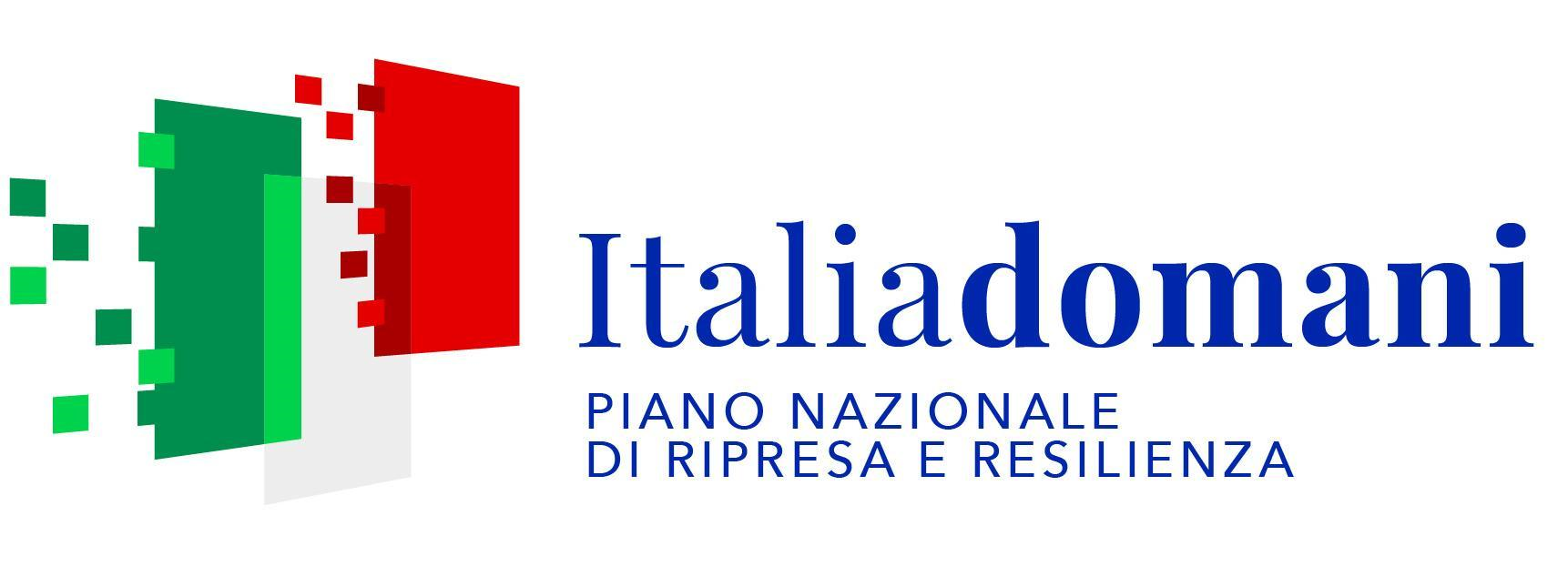}
		\hfill\includegraphics[scale=0.075]{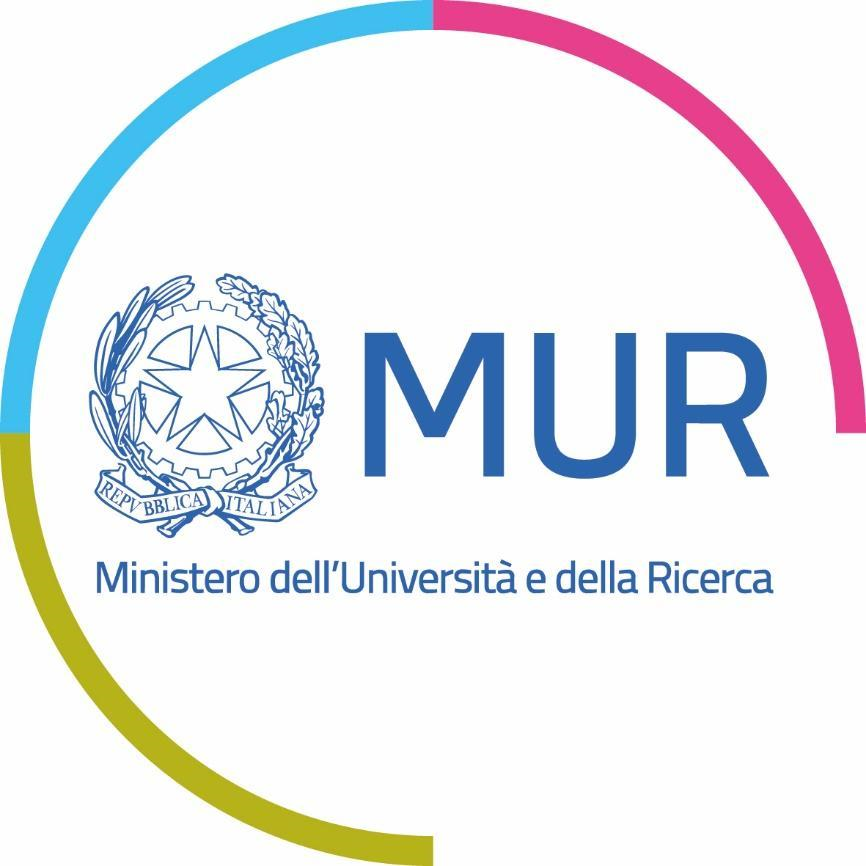}\hfill\mbox{}}
\end{figure}

\bibliographystyle{elsarticle-num}
\bibliography{Bibliography}

\begin{thebibliography}{10}
\expandafter\ifx\csname url\endcsname\relax
  \def\url#1{\texttt{#1}}\fi
\expandafter\ifx\csname urlprefix\endcsname\relax\def\urlprefix{URL }\fi
\expandafter\ifx\csname href\endcsname\relax
  \def\href#1#2{#2} \def\path#1{#1}\fi

\bibitem{Hughes2005}
T.~J.~R. Hughes, J.~A. Cottrell, Y.~Bazilevs, Isogeometric analysis: {CAD},
  finite elements, {NURBS}, exact geometry and mesh refinement, Computer
  Methods in Applied Mechanics and Engineering 194~(39) (2005) 4135--4195.
\newblock \href {https://doi.org/10.1016/j.cma.2004.10.008}
  {\path{doi:10.1016/j.cma.2004.10.008}}.

\bibitem{Evans2009}
J.~A. Evans, Y.~Bazilevs, I.~Babu\v{s}ka, T.~J.~R. Hughes, $n$-widths,
  sup-infs, and optimality ratios for the $k$-version of the isogeometic finite
  element method, Computer Methods in Applied Mechanics and Engineering 198
  (2009) 1726--1741.
\newblock \href {https://doi.org/10.1016/j.cma.2009.01.021}
  {\path{doi:10.1016/j.cma.2009.01.021}}.

\bibitem{Bressan2019}
A.~Bressan, E.~Sande, Approximation in {FEM}, {DG} and {IGA}: a theoretical
  comparison, Numerische Mathematik 143~(4) (2019) 923--942.
\newblock \href {https://doi.org/10.1007/s00211-019-01063-5}
  {\path{doi:10.1007/s00211-019-01063-5}}.

\bibitem{Takizawa2014}
K.~Takizawa, T.~Tezduyar, Space-time computation techniques with continuous
  representation in time ({ST}-{C}), Computational Mechanics 53~(1) (2014)
  91--99.
\newblock \href {https://doi.org/10.1007/s00466-013-0895-y}
  {\path{doi:10.1007/s00466-013-0895-y}}.

\bibitem{Langer2016}
U.~Langer, S.~E. Moore, M.~Neum{\"u}ller, Space-time isogeometric analysis of
  parabolic evolution problems, Computer Methods in Applied Mechanics and
  Engineering 306 (2016) 342 -- 363.
\newblock \href {https://doi.org/10.1016/j.cma.2016.03.042}
  {\path{doi:10.1016/j.cma.2016.03.042}}.

\bibitem{Montardini2018}
M.~Montardini, M.~Negri, G.~Sangalli, M.~Tani, Space-time least-squares
  isogeometric method and efficient solver for parabolic problems, Mathematics
  of Computation 89~(323) (2020) 1193--1227.
\newblock \href {https://doi.org/10.1090/mcom/3471}
  {\path{doi:10.1090/mcom/3471}}.

\bibitem{Loli2020}
G.~Loli, M.~Montardini, G.~Sangalli, M.~Tani, An efficient solver for
  space-time isogeometric {Galerkin} methods for parabolic problems, Computers
  and Mathematics with Applications 80~(11) (2020) 2586--2603.
\newblock \href {https://doi.org/10.1016/j.camwa.2020.09.014}
  {\path{doi:10.1016/j.camwa.2020.09.014}}.

\bibitem{Saade2021}
C.~Saadé, S.~Lejeunes, D.~Eyheramendy, R.~Saad, Space-{Time} {Isogeometric}
  {Analysis} for linear and non-linear elastodynamics, Computers \& Structures
  254 (2021) 106594.
\newblock \href {https://doi.org/10.1016/j.compstruc.2021.106594}
  {\path{doi:10.1016/j.compstruc.2021.106594}}.

\bibitem{Loli2023}
G.~Loli, G.~Sangalli, P.~Tesini, High-order spline upwind for space–time
  {Isogeometric Analysis}, Computer Methods in Applied Mechanics and
  Engineering 417~(11) (2023) 116408.
\newblock \href {https://doi.org/10.1016/j.cma.2023.116408}
  {\path{doi:10.1016/j.cma.2023.116408}}.

\bibitem{Brooks1982}
A.~N. Brooks, T.~J.~R. Hughes, Streamline upwind/{P}etrov-{G}alerkin
  formulations for convection dominated flows with particular emphasis on the
  incompressible {N}avier-{S}tokes equations, Computer Methods in Applied
  Mechanics and Engineering 32~(1) (1982) 199--259.
\newblock \href {https://doi.org/10.1016/0045-7825(82)90071-8}
  {\path{doi:10.1016/0045-7825(82)90071-8}}.

\bibitem{Franzone2014}
P.~Colli~Franzone, L.~F. Pavarino, S.~Scacchi, Mathematical cardiac
  electrophysiology, Vol.~13, Springer, 2014.
\newblock \href {https://doi.org/10.1007/978-3-319-04801-7}
  {\path{doi:10.1007/978-3-319-04801-7}}.

\bibitem{Quarteroni2019}
A.~Quarteroni, L.~Dedè, A.~Manzoni, C.~Vergara, Mathematical Modelling of the
  Human Cardiovascular System: Data, Numerical Approximation, Clinical
  Applications, Cambridge Monographs on Applied and Computational Mathematics,
  Cambridge University Press, 2019.
\newblock \href {https://doi.org/10.1017/9781108616096}
  {\path{doi:10.1017/9781108616096}}.

\bibitem{Franzone2012}
P.~Colli~Franzone, L.~F. Pavarino, S.~Scacchi, Mathematical and numerical
  methods for reaction-diffusion models in electrocardiology, Springer Milan,
  Milano, 2012, pp. 107--141.
\newblock \href {https://doi.org/10.1007/978-88-470-1935-5\_5}
  {\path{doi:10.1007/978-88-470-1935-5\_5}}.

\bibitem{Krishnamoorthi2013}
S.~Krishnamoorthi, M.~Sarkar, W.~S. Klug, Numerical quadrature and operator
  splitting in finite element methods for cardiac electrophysiology,
  International Journal for Numerical Methods in Biomedical Engineering 29~(11)
  (2013) 1243--1266.
\newblock \href {https://doi.org/10.1002/cnm.2573}
  {\path{doi:10.1002/cnm.2573}}.

\bibitem{Bendahmane2010}
M.~Bendahmane, R.~Bürger, R.~Ruiz-Baier, A multiresolution space-time adaptive
  scheme for the bidomain model in electrocardiology, Numerical Methods for
  Partial Differential Equations 26~(6) (2010) 1377--1404.
\newblock \href {https://doi.org/10.1002/num.20495}
  {\path{doi:10.1002/num.20495}}.

\bibitem{Franzone2006}
P.~Colli~Franzone, P.~Deuflhard, B.~Erdmann, J.~Lang, L.~F. Pavarino,
  {Adaptivity in Space and Time for Reaction-Diffusion Systems in
  Electrocardiology}, SIAM Journal on Scientific Computing 28~(3) (2006)
  942--962.
\newblock \href {https://doi.org/10.1137/050634785}
  {\path{doi:10.1137/050634785}}.

\bibitem{Southern2012}
J.~Southern, G.~Gorman, M.~Piggott, P.~Farrell, Parallel anisotropic mesh
  adaptivity with dynamic load balancing for cardiac electrophysiology, Journal
  of Computational Science 3~(1) (2012) 8--16.
\newblock \href {https://doi.org/10.1016/j.jocs.2011.11.002}
  {\path{doi:10.1016/j.jocs.2011.11.002}}.

\bibitem{Cantwell2014}
C.~D. Cantwell, S.~Yakovlev, R.~M. Kirby, N.~S. Peters, S.~J. Sherwin,
  {High-order spectral/hp element discretisation for reaction–diffusion
  problems on surfaces: Application to cardiac electrophysiology}, Journal of
  Computational Physics 257 (2014) 813--829.
\newblock \href {https://doi.org/10.1016/j.jcp.2013.10.019}
  {\path{doi:10.1016/j.jcp.2013.10.019}}.

\bibitem{Patelli2017}
A.~S. Patelli, L.~Dedè, T.~Lassila, A.~Bartezzaghi, A.~Quarteroni,
  {Isogeometric approximation of cardiac electrophysiology models on surfaces:
  An accuracy study with application to the human left atrium}, Computer
  Methods in Applied Mechanics and Engineering 317 (2017) 248--273.
\newblock \href {https://doi.org/10.1016/j.cma.2016.12.022}
  {\path{doi:10.1016/j.cma.2016.12.022}}.

\bibitem{Pegolotti2019}
L.~Pegolotti, L.~Dedè, A.~Quarteroni, {Isogeometric Analysis of the
  electrophysiology in the human heart: Numerical simulation of the bidomain
  equations on the atria}, Computer Methods in Applied Mechanics and
  Engineering 343 (2019) 52--73.
\newblock \href {https://doi.org/10.1016/j.cma.2018.08.032}
  {\path{doi:10.1016/j.cma.2018.08.032}}.

\bibitem{Rogers1994}
J.~Rogers, A.~McCulloch, A collocation-{Galerkin} finite element model of
  cardiac action potential propagation, IEEE Transactions on Biomedical
  Engineering 41~(8) (1994) 743--757.
\newblock \href {https://doi.org/10.1109/10.310090}
  {\path{doi:10.1109/10.310090}}.

\bibitem{torre2022efficient}
M.~Torre, S.~Morganti, A.~Nitti, M.~D. de~Tullio, F.~S. Pasqualini, A.~Reali,
  An efficient isogeometric collocation approach to cardiac electrophysiology,
  Computer Methods in Applied Mechanics and Engineering 393 (2022) 114782.

\bibitem{bucelli2021multipatch}
M.~Bucelli, M.~Salvador, A.~Dedè, L.~Quarteroni, Multipatch isogeometric
  analysis for electrophysiology: simulation in a human heart, Computer Methods
  in Applied Mechanics and Engineering 376 (2021) 113666.

\bibitem{Loli2022}
G.~Loli, G.~Sangalli, M.~Tani, Easy and efficient preconditioning of the
  isogeometric mass matrix, Computers \& Mathematics with Applications 116
  (2022) 245--264, new trends in Computational Methods for PDEs.
\newblock \href {https://doi.org/10.1016/j.camwa.2020.12.009}
  {\path{doi:10.1016/j.camwa.2020.12.009}}.

\bibitem{Daveiga2014}
L.~Beirão~da Veiga, A.~Buffa, G.~Sangalli, R.~Vázquez, Mathematical analysis
  of variational isogeometric methods, Acta Numerica 23 (2014) 157–287.
\newblock \href {https://doi.org/10.1017/S096249291400004X}
  {\path{doi:10.1017/S096249291400004X}}.

\bibitem{Cottrell2009}
J.~A. Cottrell, T.~J.~R. Hughes, Y.~Bazilevs, Isogeometric analysis: toward
  integration of {CAD} and {FEA}, John Wiley \& Sons, 2009.

\bibitem{Wenjian2018}
Y.~Wenjian, G.~Yu, L.~Yaohang, {Efficient Randomized Algorithms for the
  Fixed-Precision Low-Rank Matrix Approximation}, SIAM Journal on Matrix
  Analysis and Applications 39~(3) (2018) 1339--1359.
\newblock \href {https://doi.org/10.1137/17M1141977}
  {\path{doi:10.1137/17M1141977}}.

\bibitem{Vazquez2016}
R.~V{\'a}zquez, A new design for the implementation of isogeometric analysis in
  {O}ctave and {M}atlab: {G}eo{PDE}s 3.0, Computers \& Mathematics with
  Applications 72~(3) (2016) 523--554.
\newblock \href {https://doi.org/10.1016/j.camwa.2016.05.010}
  {\path{doi:10.1016/j.camwa.2016.05.010}}.

\bibitem{Sorber2014}
L.~Sorber, M.~Van~Barel, L.~De~Lathauwer,
  \href{http://www.tensorlab.net/}{Tensorlab v2.0}, Available online, January
  (2014).
\newline\urlprefix\url{http://www.tensorlab.net/}

\bibitem{Feischl2016}
M.~Feischl, G.~Gantner, A.~Haberl, D.~Praetorius, Adaptive {2D} {IGA} boundary
  element methods, Engineering Analysis with Boundary Elements 62 (2016)
  141--153.

\bibitem{Vuong2011}
A.-V. Vuong, C.~Giannelli, B.~J{\"u}ttler, B.~Simeon, A hierarchical approach
  to adaptive local refinement in isogeometric analysis, Computer Methods in
  Applied Mechanics and Engineering 200~(49-52) (2011) 3554--3567.

\bibitem{Bazilevs2007}
Y.~Bazilevs, V.~M. Calo, T.~E. Tezduyar, T.~J.~R. Hughes, Yz$\beta$
  discontinuity capturing for advection-dominated processes with application to
  arterial drug delivery, International Journal for Numerical Methods in Fluids
  54~(6-8) (2007) 593--608.

\bibitem{Charawi2017}
L.~A. Charawi,
  \href{https://www.sciencedirect.com/science/article/pii/S004578251631876X}{Isogeometric
  overlapping {S}chwarz preconditioners for the {B}idomain reaction–diffusion
  system}, Computer Methods in Applied Mechanics and Engineering 319 (2017)
  472--490.
\newblock \href {https://doi.org/https://doi.org/10.1016/j.cma.2017.03.012}
  {\path{doi:https://doi.org/10.1016/j.cma.2017.03.012}}.
\newline\urlprefix\url{https://www.sciencedirect.com/science/article/pii/S004578251631876X}

\bibitem{Willems2024}
R.~Willems, K.~L. P.~M. Janssens, P.~H.~M. Bovendeerd, C.~V. Verhoosel,
  O.~van~der Sluis, An isogeometric analysis framework for ventricular cardiac
  mechanics, Computational Mechanics 73~(3) (2024) 465--506.

\bibitem{franzone1998spread}
P.~C. Franzone, L.~Guerri, M.~Pennacchio, B.~Taccardi, Spread of excitation in
  {3}-{D} models of the anisotropic cardiac tissue. {II}. {E}ffects of fiber
  architecture and ventricular geometry, Mathematical biosciences 147~(2)
  (1998) 131--171.

\bibitem{pennacchio2011fast}
M.~Pennacchio, V.~Simoncini, Fast structured {AMG} preconditioning for the
  bidomain model in electrocardiology, SIAM Journal on Scientific Computing
  33~(2) (2011) 721--745.

\bibitem{willems2023isogeometric}
R.~Willems, E.~Kruithof, K.~L. P.~M. Janssens, M.~J.~M. Cluitmans, O.~van~der
  Sluis, P.~H.~M. Bovendeerd, C.~V. Verhoosel, Isogeometric-mechanics-driven
  electrophysiology simulations of ventricular tachycardia, in: International
  Conference on Functional Imaging and Modeling of the Heart, Springer, 2023,
  pp. 97--106.

\end{thebibliography}

\end{document}